\documentclass[11pt]{amsart} 
\usepackage[utf8]{inputenc}
\usepackage{fullpage}

\usepackage[english]{babel}
\usepackage[left=1 in, right=1 in, top=1.4 in, bottom=1.3 in]{geometry}
\usepackage{amsmath}
\usepackage{amssymb}
\usepackage{mathtools}
\usepackage[utf8]{inputenc}
\usepackage{color}
\usepackage{amsthm}
\usepackage{tikz}
\usetikzlibrary{graphs,arrows.meta}
\usepackage[shortlabels]{enumitem}
\usepackage{graphicx}
\usepackage{caption}
\usepackage{subcaption}
\usepackage{hyperref}
\usepackage{ytableau}
\usepackage{mathtools}
\usepackage[normalem]{ulem}
\usepackage{multirow} 
\usepackage{rotating}
\usepackage{array}
\usepackage{makecell}

\newcommand{\qbinom}[2]{\genfrac{[}{]}{0pt}{}{#1}{#2}_{q}}

\newcommand{\cupdot}{\mathbin{\mathaccent\cdot\cup}}
\newcommand{\vanish}[1]{}

\theoremstyle{definition}
\newtheorem{definition}{Definition}

\newtheorem{theorem}{Theorem}
\newtheorem{corollary}[theorem]{Corollary}

\newtheorem{example}{Example} 

\newcommand{\bsy}{\boldsymbol}

\newcommand{\vt}{vacillating tableau}
\newcommand{\vtx}{vacillating tableaux}
\makeatletter

\providecommand*{\bigcupdot}{%
  \mathop{%
    \vphantom{\bigcup}%
    \mathpalette\@bigcupdot{}%
  }%
}
\newcommand*{\@bigcupdot}[2]{%
  \ooalign{%
    $\m@th#1\bigcup$\cr
    \sbox0{$#1\bigcup$}%
    \dimen@=\ht0 %
    \advance\dimen@ by -\dp0 %
    \sbox0{\scalebox{2}{$\m@th#1\cdot$}}%
    \advance\dimen@ by -\ht0 %
    \dimen@=.5\dimen@
    \hidewidth\raise\dimen@\box0\hidewidth
  }%
}
\makeatother

 \usepackage[colorinlistoftodos]{todonotes}

\title{Combinatorial Identities for  Vacillating Tableaux }

\author{Zhanar Berikkyzy}
\author{Pamela E.~Harris}
\author{Anna Pun}  
\author{Catherine Yan}
\author{Chenchen Zhao} 

\address[Berikkyzy]{Department of Mathematics, Fairfield University, Fairfield, CT 06824 } 
\email{zberikkyzy@fairfield.edu}

\address[Harris]{Department of Mathematical Sciences, University of Wisconsin Milwaukee, Milwaukee, WI 53211} 
\email{peharris@uwm.edu}

\address[Pun]{Department of Mathematics, CUNY Baruch College, New York, NY 10010} 
\email{anna.pun@baruch.cuny.edu} 

\address[Yan]{Department of Mathematics, Texas A\&M University, College Station, TX 77843}
\email{huafei-yan@tamu.edu}  

\address[Zhao]{Department of Mathematics, University of Southern California, Los Angeles, CA 90089} 
\email{zhao109@usc.edu}

\date{} 

\begin{document}

\begin{abstract}\ 
Vacillating tableaux are sequences of integer partitions that satisfy specific conditions.  The concept 
of \vtx \ stems  from 
the representation theory of the partition algebra and the combinatorial theory of crossings and nestings of 
matchings and set partitions. 
In this paper, we further investigate the enumeration  of  \vtx\  
and derive multiple combinatorial identities and integer sequences  relating to the number of vacillating tableaux, simplified vacillating tableaux, and 
 limiting vacillating tableaux.  
\end{abstract}

\maketitle 

\tableofcontents

\section{Introduction}

A vacillating tableau is a sequence of integer partitions that must adhere to specific  conditions 
and can be visualized as particular walks on Young's lattice, 
 a lattice that consists of all integer partitions ordered by the containment of diagrams.
The concept of vacillating tableaux was independently introduced by two research teams around the same time.
Halverson and Lewandowski \cite{HL05} introduced one definition to provide 
combinatorial proofs of identities arising from the representation theory of the partition algebra $\mathcal{CA}_k(n)$. Meanwhile, 
 Chen, Deng, Du, Stanley, and Yan~\cite{CDDSY07} proposed another definition 
 to characterize the maximal crossings and nestings in the arc-diagrams associated with perfect matchings and set partitions of $[k]:=\{1, 2, \dots, k\}$.  
It is worth noting that these two definitions are equivalent when $n \geq 2k$. 
Additionally, in our previous work \cite{BHPYZ} 
 we introduced the concept of  \emph{limiting vacillating tableaux} 
 while studying the intricate  properties of a bijection constructed in \cite{HL05}. 
These results further motivated our investigation on the enumeration  of  \vtx, 
resulting in a number of combinatorial identities and integer sequences, which we will present here.

We begin by providing the necessary definitions.
A partition of a positive integer $n$ is a sequence  $\lambda=(\lambda_1, \dots, \lambda_t)$ of integers such that $\lambda_1 \geq \lambda_2 \geq \cdots \geq \lambda_t >0$ and $|\lambda|\coloneqq\lambda_1+ \cdots +\lambda_t =n$. 
We also say that the size of $\lambda$ is $n$ and denote it as  $\lambda \vdash n$. 
In addition, we let the empty partition $\emptyset$ to be the only integer partition of $0$.

A partition $\lambda$ is 
visually represented by the Young diagram, which contains $\lambda_j$ boxes in the $j$-th row.  We adopt the English notation in which the diagrams are aligned in the upper-left corner.
 A \emph{Young tableau} of shape $\lambda$ is an array obtained by filling each box of the Young diagram of $\lambda$ with an integer. 
 A Young tableau is \emph{semistandard} if the entries are  weakly increasing in every row and strictly increasing in every column. 
 A semistandard Young tableau (SSYT) is \emph{partial} if the entries are all distinct, and we call it a partial tableau. 
 The \emph{content} of a Young tableau $T$, denoted as content$(T)$,  is the multiset of all the entries in $T$.
If the content of a SSYT $T$ with shape $\lambda \vdash n$ is exactly $[n]$, then we say $T$ is 
  a \emph{standard Young tableau} (SYT).  
  Throughout this paper, we use $f^\lambda$ to denote the number of SYTs of shape $\lambda$. By convention, we set $f^\emptyset = 1$.
  
Below is the definition introduced in \cite{HL05}. To emphasize the dependency on the parameter $n$, we call such tableaux $n$-vacillating tableaux. 

 \begin{definition} \cite{HL05}
For integers $k \geq 0$ and $n \geq 1$, an \emph{$n$-vacillating tableau of shape $\lambda$ and length $2k$} is a sequence of $2k+1$ integer partitions 
\[
((n)=\lambda^{(0)}, \lambda^{(\frac{1}{2})}, \lambda^{(1)}, \lambda^{(1\frac{1}{2})}, \dots, \lambda^{(k-\frac{1}{2})}, \lambda^{(k)}=\lambda) 
\]
so that for each $j = 0,1, \dots,  k -1$, 
\begin{enumerate}[(a)]
    \item $\lambda^{(j)} \supseteq  \lambda^{(j+\frac12)}$ and 
    $| \lambda^{(j)} / \lambda^{(j+\frac12)}|=1$, 
    \item $\lambda^{(j+\frac12)}\subseteq \lambda^{(j+1)}$ and 
    $|\lambda^{(j+1)}/ \lambda^{(j+\frac12)}|=1$. 
    \end{enumerate} 
    \end{definition} 
In other words, an $n$-\vt\ of shape $\lambda$ is a walk on Young's lattice from $(n)$ to $\lambda$, where a box is removed in each odd step and  added in each even step.

Let 
$
\Lambda_n^k = \left\{  \lambda \vdash n:  \   \lambda_1 \geq n -k\right\}
$
and $\mathcal{VT}_{n,k}(\lambda)$ be the set of all $n$-vacillating tableaux of shape $\lambda$ and  length $2k$. 
Note that for any \vt \ in $\mathcal{VT}_{n,k}(\lambda)$,  $\lambda^{(j)} \in \Lambda_{n}^k$ and $\lambda^{(j+\frac12)} \in \Lambda_{n-1}^k$.
Let $m_{n,k}^{\lambda}$ be the cardinality of $\mathcal{VT}_{n,k}(\lambda)$.  When $n \geq 2k$, the value of $m_{n,k}^{\lambda}$ does not depend on $n$. In that case, 
there is an equivalent notion of vacillating tableau that does not use the parameter $n$, which was 
mentioned in \cite{HL05} but originally introduced in \cite{CDDSY07} to study maximal monotone 
substructures in matchings and set partitions. 
To distinguish from the $n$-vacillating tableau, we call this second notion the \emph{simplified vacillating tableau}. 
Given an integer partition $\lambda=(\lambda_1, \lambda_2, \dots, \lambda_t) \vdash n$, 
let $\lambda^*$ be obtained from $\lambda$ by removing its first part $\lambda_1$, i.e., 
$\lambda^*=(\lambda_2, \dots, \lambda_t) \vdash (n-\lambda_1)$.
For an $n$-vacillating tableau $( \lambda^{(j)}: j=0, \frac{1}{2}, 1, 1\frac{1}{2}, \dots, k)$ in $\mathcal{VT}_{n,k}(\lambda)$, the corresponding \emph{simplified vacillating tableau}  is the sequence 
$( \mu^{(j)}: \mu^{(j)}= (\lambda^{(j)})^* \text{ for  } j=0, \frac{1}{2}, 1, 1\frac{1}{2}, \dots, k)$. 
One can also define the simplified vacillating tableau  directly in terms of integer partitions. 

\begin{definition} \cite{CDDSY07}  \label{def:svt} 
A \emph{simplified vacillating tableau}  of shape $\mu$  and length $2k$  is a 
sequence  of $2k+1$ integer partitions 
$$
(\mu^{(0)} = \emptyset,  \mu^{(\frac{1}{2})}, \mu^{(1)}, \mu^{(1\frac{1}{2})}, \dots, \mu^{(k-\frac{1}{2})},
\mu^{(k)} = \mu) 
$$ 
so that for each integer $j = 0, 1, \dots, k-1$, 
\begin{enumerate}[(a)]
\item $\mu^{(j)} \supseteq  \mu^{(j+\frac12)}$ and 
    $| \mu^{(j)} / \mu^{(j+\frac12)}|=0$ or $1$, 
    \item $\mu^{(j+\frac12)}\subseteq \mu^{(j+1)}$ and 
    $|\mu^{(j+1)}/ \mu^{(j+\frac12)}|=0$ or $1$. 
\end{enumerate}
 \end{definition} 

Note that the definition forces $\mu^{(\frac{1}{2})}=\emptyset$ and $|\mu| \leq k$. 
Equivalently, a simplified \vt\  of shape $\lambda$ is a walk on Young's lattice from $\emptyset$ to $\lambda$ where in each odd step either zero or one box is removed, and in each even step either zero or one box is added.
Let $\mathcal{SVT}_k(\mu)$ be the set of all
simplified vacillating tableaux of shape $\mu $ and length $2k$, and let $g_k(\mu) $  be the cardinality of $\mathcal{SVT}_k(\mu)$.  
In general,   $m_{n,k}^\lambda \leq g_k(\lambda^*) $,   
with equality holding when $n \geq 2k$.

\begin{figure}[ht] 
\begin{center} 
    \begin{tikzpicture}
     \node at (0,0) {$k=0:$};  
     \node at (2,0) {$\emptyset$}; 
     \node[red] at (2.3,-0.2) {\small $1$}; 
     \draw[->] (2, -0.3)--(2, -1.1); 
     
     \node at (0, -1.5) {$k=\frac{1}{2}:$}; 
     \node at (2, -1.5) {$\emptyset$}; 
     \node[red] at (2.3, -1.7) {\small{$1$}}; 
     \draw[->] (2, -1.8)--(2, -2.6); 
     \draw[->] (2.2, -1.8)--(3.2, -2.6); 
     
     \node at (0, -3) {$k=1$:}; 
     \node at (2, -3) {$\emptyset$};
     \draw (3.0, -3.2) rectangle (3.4, -2.8); 
     \node[red] at (2.3, -3.2) {\small $1$}; 
     \node[red] at (3.7, -3.2) {\small$1$}; 
     \draw[->] (2,-3.3)--(2,-4.1); \draw[->](3.2, -3.3)--(3.2, -4.1); 
     \draw[->] (3, -3.3)--(2.2,-4.1); 
     
     \node at (0, -4.5) {$k=1\frac{1}{2}:$}; 
     \node at (2, -4.5) {$\emptyset$};
     \node[red] at (2.3, -4.7) {\small $2$}; 
     \draw (3.0, -4.7) rectangle (3.4, -4.3); 
     \node[red] at (3.7, -4.7) {\small$1$}; 
     \draw[->] (2, -4.8)--(2, -5.6); \draw[->] (2.2, -4.8)--(3, -5.6); 
     \draw[->] (3.2, -4.8)--(3.2, -5.6);
     \draw[->] (3.4, -4.8)--(4.8, -5.6); 
     \draw[->] (3.6, -4.8) --(7, -5.6);

     \node at (0, -6) {$k=2:$}; 
     \node at (2, -6) {$\emptyset$};
     \node[red] at (2.2, -6.2) {\small$2$}; 
     \draw (3.2, -6.2) rectangle (3.6, -5.8); 
     \node[red] at (3.8, -6.2) {\small$3$}; 
     \draw (4.8, -6.2) rectangle (5.6, -5.8); 
     \draw (5.2, -5.8)--(5.2, -6.2); 
     \node[red] at (5.8, -6.2) {\small$1$}; 
     \draw (7, -5.8) rectangle (7.4, -6.6); 
     \draw (7, -6.2)--(7.4, -6.2); 
     \node[red] at (7.6, -6.6) {\small$1$}; 
     \draw[->] (2, -6.5)--(2, -7.3); 
     \draw[->] (3.2, -6.5)--(2.2, -7.3); 
     \draw[->] (3.4, -6.5)--(3.4, -7.3); 
     \draw[->] (5.2, -6.5)--(3.6, -7.3); 
     \draw[->] (5.4, -6.5)--(5.4, -7.4); 
     \draw[->] (7, -6.7)--(5.8, -7.3); 
     \draw[->] (7.2, -6.7)--(7.2,-7.3);

     \node at (0, -7.9) {$k=2\frac{1}{2}:$}; 
     \node at (2, -7.9) {$\emptyset$}; 
     \node[red] at (2.3, -8.2) {\small$5$}; 
     \draw (3.2, -8.1) rectangle (3.6, -7.7); 
     \node[red] at  (3.9, -8.2) {\small$5$}; 
     \draw (4.8, -8.1) rectangle (5.6, -7.7); 
     \draw (5.2, -8.1)--(5.2, -7.7); 
      \node[red] at (5.9, -8.2) {\small$1$}; 
     \draw (7, -7.6) rectangle (7.4, -8.4); 
       \draw (7, -8)--(7.4, -8); 
       \node[red] at (7.7, -8.2) {\small$1$}; 
      \draw[->] (2, -8.3)--(2, -9.1);   \draw[->] (2.2,-8.3)--(3.2, -9.1); 
      \draw[->] (3.4, -8.3)--(3.4,-9.1);  \draw[->] (3.6,-8.3)--(4.8, -9.1); \draw[->] (3.6, -8.3)--(7,-9.1); 
      \draw[->] (5.4, -8.3)--(5.4, -9.1); 
      \draw[->] (5.6, -8.3)--(8.6, -9.2); \draw[->] (5.8, -8.3)--(10.5, -9.2); 
      \draw[->] (7.2, -8.4) -- (7.2,-9.1); 
      \draw[->] (7.9, -8.3)--(11, -9.1); 
      \draw[->] (8.1, -8.3)--(12.6, -9.1);

     \node at (0, -9.6) {$k=3:$}; 
     \node at (2, -9.6) {$\emptyset$}; 
       \node[red] at (2.2,-9.8) {\small$5$}; 
     \draw (3.2, -9.8) rectangle (3.6, -9.4); 
      \node[red] at (3.8, -9.8) {\small$10$}; 
     \draw (4.8, -9.8) rectangle (5.6, -9.4); 
      \draw (5.2, -9.8)--(5.2, -9.4); 
      \node[red] at (5.8, -9.8) {\small$6$}; 
     \draw (7, -9.4) rectangle (7.4, -10.2); 
       \draw (7,-9.8)--(7.4, -9.8); 
       \node[red] at (7.6, -10.2) {\small$6$}; 
     \draw (8.8, -9.8) rectangle (10, -9.4); 
        \draw (9.2,-9.8)--(9.2,-9.4) (9.6,-9.8)--(9.6,-9.4); 
        \node[red] at (10.3, -9.8) {\small$1$}; 
    \draw (10.8, -9.8) rectangle (11.6,-9.4); 
    \draw (10.8,-9.4) rectangle (11.2,-10.2); 
    \node[red] at (11.6, -10.2) {\small$2$}; 
    \draw (12.4, -9.4) rectangle (12.8, -10.6); 
      \draw (12.4, -9.8) rectangle (12.8, -10.2); 
      \node[red] at (13.1, -10.6) {\small$1$}; 
    \end{tikzpicture}
    \caption{Simplified vacillating tableaux of length $2k$, up to $k=3$.}  \label{fig:svt} 
    \end{center} 
    \end{figure}
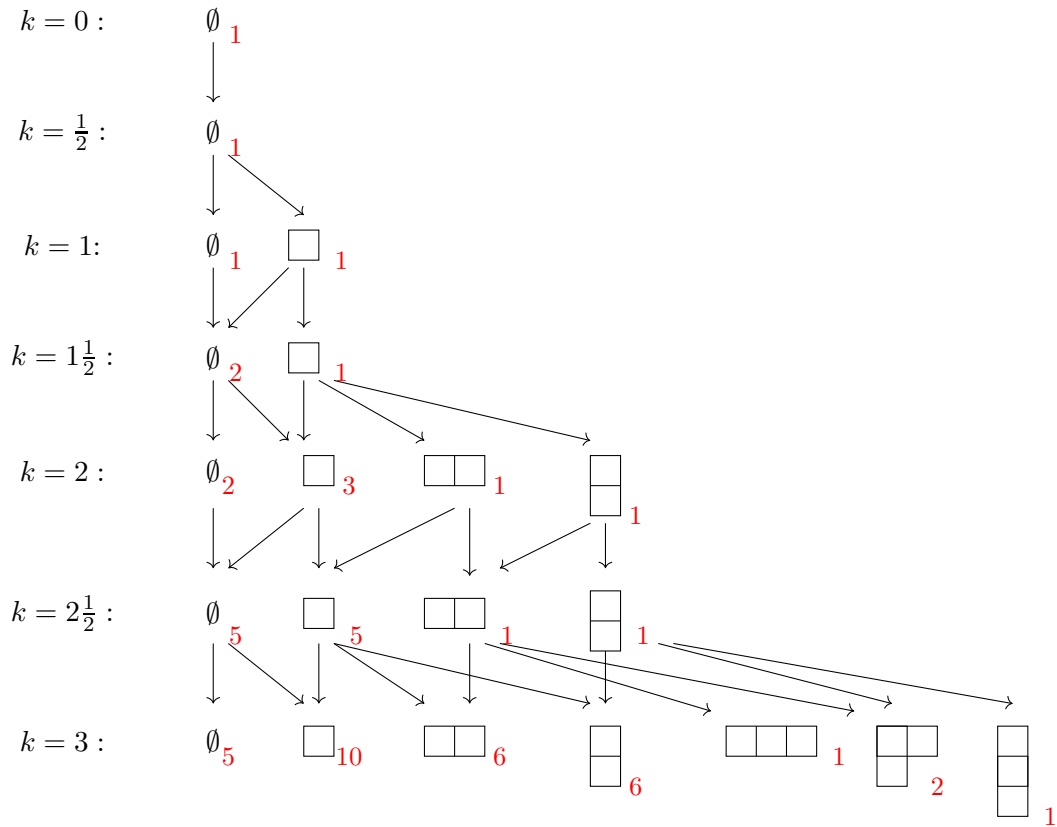

The diagram in Figure~\ref{fig:svt}  gives the simplified vacillating tableaux of length $2k$, for $0 \leq  k \leq 3$.  The number next to each integer partition $\mu$ is the value of $g_k(\mu)$.

A limiting vacillating tableau is a special kind of simplified vacillating tableau tableau, which was  introduced in \cite{BHPYZ}  from the study of images of a bijection from \cite{HL05} between  sequences in $\{ (i_1, \dots, i_k):\  1 \leq i_j \leq n\}$ and pairs $(T, P)$,  where $T$ is a SYT of 
some shape $\lambda \in \Lambda_n^k$,  and $P$ is an $n$-\vt \ of shape $\lambda$ and length $2k$.

\begin{definition} \cite{BHPYZ} 
A \emph{limiting vacillating tableau}  of shape $\mu$  and length $2k$  is a 
sequence  of $2k+1$ integer partitions 
$$
(\mu^{(0)} = \emptyset,  \mu^{(\frac{1}{2})}, \mu^{(1)}, \mu^{(1\frac{1}{2})}, \dots, \mu^{(k-\frac{1}{2})},
\mu^{(k)} = \mu) 
$$ 
so that for each integer $j=0, 1, \dots, k-1$, 
\begin{enumerate}[(a)]
\item $\mu^{(j)} \supseteq  \mu^{(j+\frac12)}$ and 
    $| \mu^{(j)} / \mu^{(j+\frac12)}|=0$ or $1$, 
    \item $\mu^{(j+\frac12)}\subseteq \mu^{(j+1)}$ and 
    $|\mu^{(j+1)}/ \mu^{(j+\frac12)}|=1$. 
     \label{b-part}
\end{enumerate}
\end{definition}

Equivalently, a limiting \vt\  of shape $\lambda$ is a walk on Young's lattice from $\emptyset$ to $\lambda$ where in each odd step either zero or one box is removed, and in each even step exactly one box is added.

\begin{figure}[ht] 
\begin{center} 
    \begin{tikzpicture}
     \node at (0,0) {$k=0:$};  
     \node at (2,0) {$\emptyset$}; 
     \node[red] at (2.3,-0.2) {\small $1$}; 
     \draw[->] (2, -0.3)--(2, -1.1); 
     
     \node at (0, -1.5) {$k=\frac{1}{2}:$}; 
     \node at (2, -1.5) {$\emptyset$}; 
     \node[red] at (2.3, -1.7) {\small{$1$}}; 
     \draw[->] (2, -1.8)--(2, -2.6); 
     
     \node at (0, -3) {$k=1$:}; 
     \draw (1.8, -3.2) rectangle (2.2, -2.8); 
     \node[red] at (2.5, -3.2) {\small $1$}; 
     \draw[->] (2,-3.3)--(2,-4.1); \draw[->](2.2, -3.3)--(3.2, -4.1); 
     
     \node at (0, -4.5) {$k=1\frac{1}{2}:$}; 
     \node at (2, -4.5) {$\emptyset$};
     \node[red] at (2.3, -4.7) {\small $1$}; 
     \draw (3.0, -4.7) rectangle (3.4, -4.3); 
     \node[red] at (3.7, -4.7) {\small$1$}; 
     \draw[->] (2, -4.8)--(2, -5.6); 
     \draw[->] (3.2, -4.8)--(3.2, -5.6);
     \draw[->] (3.4, -4.8)--(4.8, -5.6);

     \node at (0, -6) {$k=2:$}; 
     \draw (1.8, -6.2) rectangle (2.2, -5.8); 
     \node[red] at (2.5, -6.2) {\small$1$}; 
     \draw (3.0, -6.2) rectangle (3.8, -5.8); 
     \draw (3.4, -5.8)--(3.4, -6.2); 
     \node[red] at (4.1, -6.2) {\small$1$}; 
     \draw (4.6, -5.8) rectangle (5.0, -6.6); 
     \draw (4.6, -6.2)--(5.0, -6.2); 
     \node[red] at (5.3, -6.6) {\small$1$}; 
     \draw[->] (2, -6.5)--(2, -7.3); 
     \draw[->] (2.2, -6.5)--(3.2, -7.3); 
     \draw[->] (3.4, -6.5)--(3.4, -7.3); 
     \draw[->] (3.6, -6.5)--(5.2, -7.3); 
     \draw[->] (5.0, -6.7)--(7.2, -7.4); 
     \draw[->] (4.8, -6.7)--(3.6, -7.3);

     \node at (0, -7.9) {$k=2\frac{1}{2}:$}; 
     \node at (2, -7.9) {$\emptyset$}; 
     \node[red] at (2.3, -8.2) {\small$1$}; 
     \draw (3.2, -8.1) rectangle (3.6, -7.7); 
     \node[red] at  (3.9, -8.2) {\small$3$}; 
     \draw (4.8, -8.1) rectangle (5.6, -7.7); 
     \draw (5.2, -8.1)--(5.2, -7.7); 
      \node[red] at (5.9, -8.2) {\small$1$}; 
     \draw (7, -7.6) rectangle (7.4, -8.4); 
       \draw (7, -8)--(7.4, -8); 
       \node[red] at (7.7, -8.2) {\small$1$}; 
      \draw[->] (2, -8.3)--(2, -9.1);  
      \draw[->] (3.4, -8.3)--(3.4,-9.1); 
      \draw[->] (3.6,-8.3)--(4.6, -9.1); 
      \draw[->] (5.2, -8.3)--(5.8, -9.1); 
      \draw[->] (5.4, -8.3)--(7.6, -9.1); 
      \draw[->] (7.2, -8.5)--(7.8, -9.1); 
      \draw[->] (7.4, -8.5)--(9.6, -9.1);

     \node at (0, -9.6) {$k=3:$}; 
     \draw (1.8, -9.8) rectangle (2.2, -9.4); 
      \node[red] at (2.5, -9.8) {\small$1$}; 
     \draw (3.0, -9.8) rectangle (3.8, -9.4); 
      \draw (3.4, -9.8)--(3.4, -9.4); 
      \node[red] at (4.1, -9.8) {\small$3$}; 
     \draw (4.6, -9.4) rectangle (5.0, -10.2); 
       \draw (4.6,-9.8)--(5.0, -9.8); 
       \node[red] at (5.3, -10.2) {\small$3$}; 
     \draw (5.8, -9.8) rectangle (7, -9.4); 
        \draw (6.2,-9.8)--(6.2,-9.4) (6.6,-9.8)--(6.6,-9.4); 
        \node[red] at (7.3, -9.8) {\small$1$}; 
    \draw (7.8, -9.8) rectangle (8.6,-9.4); 
    \draw (7.8,-9.4) rectangle (8.2,-10.2); 
    \node[red] at (8.6, -10.2) {\small$2$}; 
    \draw (9.4, -9.4) rectangle (9.8, -10.6); 
      \draw (9.4, -9.8) rectangle (9.8, -10.2); 
      \node[red] at (10.1, -10.6) {\small$1$}; 
    \end{tikzpicture}
    \caption{Limiting vacillating tableaux of length $2k$, up to $k=3$.}  \label{fig:lvt} 
    \end{center} 
    \end{figure}
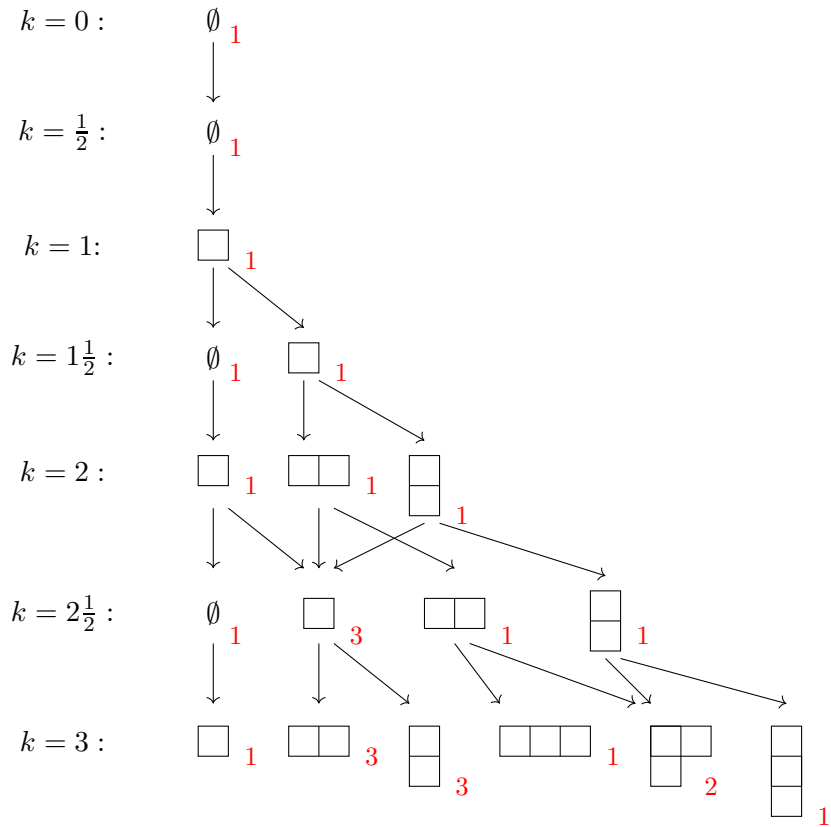

Figures~\ref{fig:lvt} gives the limiting vacillating tableaux of length $2k$, for $0 \leq  k \leq 3$.  The number next to each integer partition $\mu$ is the number of limiting vacillating tableaux  ending at $\mu$,  which is denoted by $a_k(\mu)$. 


The above definitions  suggest that we can also consider the numbers of vacillating tableaux and their simplified analogs of odd length, 
namely, $m_{n, k+\frac{1}{2}}^\lambda$ and $g_{k+\frac{1}{2}}(\mu)$, where 
$g_{k+\frac{1}{2}}(\mu)=m_{n, k+\frac{1}{2}}^\lambda$ if $n\geq 2k$ and $\mu=\lambda^*$. Similarly, we define 
$a_{k+\frac{1}{2}}(\mu)$ as the number of limiting \vt x of shape $\mu$ and length $2k+1$.

In Section \ref{sec:bkgd},  we review some combinatorial algorithms and bijections related to \vtx,  as well as identities 
of $m_{n,k}^\lambda$ arising from the representation theory of the partition algebra. 
Motivated by these results, we  consider various summations concerning $g_k(\mu), g_{k+\frac{1}{2}}(\mu),$ $a_k(\mu)$ and $a_{k+\frac{1}{2}}(\mu)$, and  obtain new identities and 
present combinatorial interpretations. Sections \ref{sec:svt1}  and  \ref{sec:svt2} focus on simplified \vtx. We will specify when there is an analogous result for $m_{n,k}^\lambda$ for general $n, k$, (without the constraint that $n \geq 2k$).  Sections \ref{sec:lvt1} and \ref{sec:lvt2} focus on limiting \vtx. 
More precisely, we consider the following integer sequences. 
\begin{itemize} 
 \item The number of  simplified vacillating tableaux of lengths $2k$  and $2k+1$: 
 $g_k=\sum_{\mu} g_k(\mu)$    and 
  $g_{k+\frac{1}{2}}=\sum_\mu g_{k+\frac{1}{2}}(\mu)$. 
\item  The number of  limiting vacillating tableaux of lengths $2k$ and $2k+1$: 
$a_k=\sum_{\mu} a_k(\mu)$   and 
  $a_{k+\frac{1}{2}}=\sum_\mu a_{k+\frac{1}{2}}(\mu)$.
  \item Sums of products of numbers of  simplified \vtx \  of shape $\mu$ and SYTs of shape $\mu$: 
   $u_k=\sum_\mu g_k(\mu)f^\mu$  and  $u_{k+\frac{1}{2}}=\sum_\mu g_{k+\frac{1}{2}}(\mu) f^\mu$. 
   \item Sums of products of numbers of limiting \vtx  \ of shape $\mu$ and SYTs of shape $\mu$: 
   $v_k=\sum_\mu a_k(\mu) f^\mu$ and $v_{k+\frac{1}{2}}=\sum_\mu a_{k+\frac{1}{2}}(\mu) f^\mu$. 
\end{itemize} 
In addition, we present results for the products of 
the numbers of \vtx\ 
with each other, with $f^\mu$, and with Schur function $s_\mu$.      
In Section \ref{sec:final}, we discuss the relation between the aforementioned integer sequences and conclude this paper  with some final remarks and future research projects.

\section{Algorithms and Bijections on Tableaux } \label{sec:bkgd} 

In this section, we review combinatorial algorithms and bijections related to \vtx,  as well as identities 
of $m_{n,k}^\lambda$ arising from  representation theory. 

\subsection{The RSK insertion algorithm  } 
A key ingredient  needed in defining the bijections  is the RSK row insertion algorithm,  which constructs a pair of tableaux of the same shape from an integer sequence, or more generally, a two-line array of integers. There are two major components: the row insertion procedure, and the RSK algorithm that is an iteration of row insertions.

\medskip

\noindent 
\underline{The RSK row insertion}. Let $T$ be a SSYT  of  shape $\lambda$ and $x$ be an integer. The operation $x \xrightarrow{RSK}  T$ is defined as follows. 
\begin{enumerate}[(a)]
    \item Let $R$ be the first row of $T$. 
    \item While $x$ is less than some entries in $R$, do 
       \begin{enumerate}[label=\roman*)]
           \item Let $y$ be the smallest entry of $R$ greater than $x$; 
           \item Replace $y \in R$ with $x$;
           \item Let $x \coloneqq y$ and let $R$ be the next row. 
       \end{enumerate}
    \item Place $x$ at the end of $R$ (which is possibly empty). 
\end{enumerate}
The result is a semistandard  tableau of shape $\mu$ such that $|\mu /\lambda|=1$. For each occurrence of step (b), we say that $x$ \emph{bumps} $y$ to the next row.  

We will use Knuth's construction \cite{Knuth70} that iterates the above  insertion procedure on a two-line 
array of integers  
\begin{equation}  \label{2-line-array}
\left( 
\begin{array}{cccc}
    u_1 & u_2 & \cdots & u_n \\ 
    v_1  & v_2 & \cdots  & v_n 
 \end{array}
\right),
\end{equation}
where $(u_j, v_j)$  are arranged in non-decreasing lexicographic order from left to right, that is, $u_1 \leq u_2 \leq \cdots \leq u_n$ and $v_j\leq v_{j+1}$ if $u_j=u_{j+1}$. 
\medskip 

\noindent 
\underline{The RSK algorithm.} \cite{Knuth70} \  
Given the two-line array \eqref{2-line-array}, construct a pair of Young tableaux $(P,Q)$ of the same shape by 
starting with $P=Q=\emptyset$. For $j=1, 2, \dots, n$,
\begin{enumerate}[(a)]
    \item Insert $v_j$ into tableau $P$ using the row insertion procedure. This 
    operation adds a new box to the shape of $P$. Assume the new box is 
    at the end of the $i$-th row of $P$.
    \item Add  a new box with entry $u_j$ at the  end of the $i$-th row  of $Q$.  
\end{enumerate}
We call $P$ the \emph{insertion tableau} and $Q$ the \emph{recording tableau}. 

\medskip
Let $A$ and $B$ be two totally ordered alphabets.   If $u_j \in A$ and $v_j \in B$ for all $j$, we say that the two-line array \eqref{2-line-array} is a generalized permutation from $A$ to $B$.  

\begin{theorem}[Knuth] 
    There is a one-to-one correspondence between generalized permutations from $A$ to $B$ and pairs of SSYTs $(P,Q)$ of the same shape, where content$(P) \subseteq B$ and content$(Q) \subseteq A$.    
\end{theorem}

Knuth's correspondence, when restricted to 
two special families of two-line arrays,  gives the correspondence discovered by 
Robinson \cite{Robinson38} and Schensted \cite{Sch61}.  
  
\begin{enumerate}[(i)]
\item When $(u_1\  u_2\ \dots\ u_n)=(1\  2\ \dots\ n)$  and $(v_1\  \dots\  v_n)$ ranges over all permutations of $[n]$,  the correspondence gives a bijection between permutations of length $n$ and pairs of SYTs of the same shape. 
\item 
When $(u_1\  u_2 \ \dots\ u_n)=(1\  2\  \dots\ n)$ and $v_i \in \mathbb{Z}^+$,  the correspondence gives a bijection between integer sequences of length $n$
and pairs of Young tableaux of the same shape 
$\lambda \vdash n$, where $P$ is a SSYT  with content in $\mathbb{Z}^+$, and $Q$ is a SYT. 
\end{enumerate}

\subsection{The delete-insert process} 
The notion of $n$-\vt\  was introduced by Halverson and Lewandowski \cite{HL05} as part of the image of a delete-insert process,   which gave a combinatorial proof of Identity \eqref{n-vt} that reflects the Schur-Weyl duality between the symmetric group algebra and the partition algebra. 
\begin{theorem} \cite{HL05}
For all integers $n \geq 1$ and $k \geq 0$, 
 \begin{equation} \label{n-vt} 
        n^k = \sum_{\lambda \in \Lambda_n^k} f^\lambda m_{n,k}^\lambda,  
    \end{equation}
    where $f^\lambda$ is the number of SYTs  of shape $\lambda$ and 
    $m_{n,k}^\lambda$ is the number of $n$-\vtx\ of shape $\lambda $ and length $2k$. 
\end{theorem}

Halverson and Lewandowski's proof  uses a combination of the RSK
row insertion algorithm and jeu de taquin, which removes a box from a given Young tableau, to  construct a bijection  $DI_n^k$ from integer sequences in $[n]^k$ to the disjoint union 
$\cupdot_{\lambda \in \Lambda_n^k} \mathcal{SYT}(\lambda) \times \mathcal{VT}_{n,k}(\lambda)$, where $\mathcal{SYT}(\lambda)$ is the set of SYTs of shape $\lambda$.  
We describe jeu de taquin and the construction of $DI_n^k$ below. 
\medskip 

\noindent 
\underline{Jeu de taquin}. \ 
Let $T$ be a partial tableau of  shape $\lambda$.  
Let $x$ be an entry in $T$.  The following operation will delete $x$ from $T$ and yield a partial tableau $S$.
\begin{enumerate}[(a)]
    \item Let $c=T_{i,j}$  be the box of $T$ containing $x$, 
    which is  the $j$-th box in the $i$-th row of $T$. 
    \item While either  $T_{i+1, j}$ or $T_{i,j+1}$ exists, 
    do the following steps. 
      \begin{enumerate}[label=\roman*)]
          \item Let $c'$ be the box containing $\min\{T_{i+1,j}, T_{i, j+1}\}$; if only one of $T_{i+1,j}, T_{i, j+1}$ exists, then the minimum is taken to be that single entry; 
          \item Exchange the positions of $c$ and $c'$. 
      \end{enumerate}
    \item Delete $c$. 
\end{enumerate}
We denote this process by $x \xleftarrow{ jdt } T$. 

\medskip

\noindent \underline{The bijection $DI_n^k$.}\ 
For any integer sequence $\bsy{i}=(i_1, i_2, \dots, i_k) \in [n]^k$, we define a sequence of  tableaux recursively: 
the $0$-th tableau $T^{(0)}$ is the unique SYT of shape $(n)$,   namely, 
\begin{center} 
\begin{tikzpicture}
\node[left] at (0,0) {$T^{(0)} =$}; 
\draw (0,-.25) rectangle (3,.25); 
\draw (0.5, -.25)--(0.5,0.25) (1,-0.25)--(1,0.25) (2.5, -0.25)--(2.5, 0.25); 
\node at (0.25,0) {$1$}; 
\node at (0.75, 0) {$2$}; 
\node at (1.75,0) {$\dots$}; 
\node at (2.75,0) {$n$}; 
\node at (3.2, -0.2) {.}; 
\end{tikzpicture}
\end{center} 
For integers $j=0, 1, \dots, k-1$, the partial tableaux $T^{(j+\frac{1}{2})}$ and $T^{(j+1)}$ are defined by
\begin{eqnarray}
T^{(j+\frac{1}{2})} &= &\left( i_{j+1} \xleftarrow{jdt} T^{(j)} \right), \\ T^{(j+1)}  & = & \left( i_{j+1} \xrightarrow{RSK} T^{(j+\frac{1}{2})} \right). 
\end{eqnarray}
It follows that for each integer index $j$, $T^{(j)}$ is a SYT of shape $\lambda^{(j)} \in \Lambda_n^k$, and  $T^{(j+\frac{1}{2})}$ is a partial 
tableau of shape  $\lambda^{(j+\frac{1}{2})} \in \Lambda_{n-1}^k$.   The map $DI_n^k$ is defined by
\begin{equation}
  \bsy{i}=  (i_1, i_2, \dots, i_k) \xrightarrow{DI_n^k} (T_\lambda, P_\lambda), 
\end{equation}
where
$\lambda= \lambda^{(k)} \in \Lambda_n^k$, 
$T_\lambda= T^{(k)} \in \mathcal{SYT}(\lambda)$ 
and $P_\lambda=(\lambda^{(0)}, \lambda^{(\frac{1}{2})}, \lambda^{(1)}, 
\lambda^{(1\frac{1}{2})}, \dots, \lambda^{(k)})$ is an $n$-vacillating  tableau in $\mathcal{VT}_{n,k}(\lambda)$.

Note that  $DI_n^k$ is well-defined without the condition  $n \geq 2k$. It is a bijection since both the row insertion and 
the jeu de taquin are invertible. 
 We prove that Identity \eqref{n-vt} can be extended to \vtx\  of odd length. 

\begin{corollary} For all integers $n \geq 1$ and $k \geq 0$, 
\begin{equation}
   n^k=  \sum_{\lambda \in \Lambda_{n-1}^k} f^\lambda m_{n,k+\frac{1}{2}}^\lambda . 
 \end{equation}   
\end{corollary}
\begin{proof}  
 Applying the delete-insert process 
 to the integer 
  sequence $(i_1, i_2, \dots, i_k, n)$, we get a 
  sequence of  tableaux $\{T^{(r)}: r=0, \frac{1}{2}, 1, 1\frac{1}{2}, 2, \dots, k+\frac{1}{2}, k+1\}$. Since the last entry of the integer sequence is $n$, $T^{(k+1)}$ is obtained from $T^{(k +\frac{1}{2})}$ by adding $n$ at the end of the first row. 
  Therefore, we can ignore the last tableau $T^{(k+1)}$, and the image of $DI_n^{k+1}$
  on $(i_1, i_2, \dots, i_k, n)$ is completely determined by $T^{(k+\frac{1}{2})}$, a SYT of shape $\lambda^{(k+\frac{1}{2})} \in \Lambda_{n-1}^k$,  and the sequence of shapes of $(T^{(0)}, T^{(\frac{1}{2})}, T^{(1)}, T^{(1\frac{1}{2})}, \dots, T^{(k+\frac{1}{2})})$, which is an $n$-vacillating tableau of shape $\lambda^{(k+\frac{1}{2})}$ and length $2k+1$.   Hence the map $DI_{n}^{k+1}$ induces a bijection between the integer sequences of the form $(i_1, i_2, \dots, i_k, n)$ and the disjoint union 
  \[
   \bigcupdot_{\lambda \in \Lambda_{n-1}^k}
    \mathcal{SYT}(\lambda) \times \mathcal{VT}_{n, k+\frac{1}{2}}(\lambda),  
  \]
  which proves the corollary. 
\end{proof}

Recently Krattenthaler \cite{Kratten23}  proved a generalization of \eqref{n-vt} using growth diagrams: 
\begin{equation}
    n^k =\sum_{\lambda \vdash n} f^\lambda 
    m_{\mu}^\lambda(k), 
\end{equation}
where  $\mu$ is any fixed integer partition of $n$ and $m_\mu^\lambda(k)$ is the number of $n$-vacillating tableaux from shape $\mu$ to shape $\lambda$ in $2k$ steps.  In this paper, we only consider $n$-\vtx \ starting at the 
shape $(n)$ and simplified \vtx \ starting at the empty shape.

\subsection{Vacillating tableaux and set partitions}  \label{ssec:vt} \mbox{}  

While the map $DI_n^k$ is the main tool to study the combinatorial properties of $n$-\vtx, another bijection  
based on row insertion works better with simplified \vtx.   
It is the bijection $\psi$ defined in \cite[Section 2]{CDDSY07}, which extends  a map in  
\cite[Lemma 2.2]{Sunderam90} on up-down tableaux of arbitrary shapes. 
The bijection $\psi$ maps  simplified \vtx \ of shape $\mu$ and length $2k$
to pairs 
$(\bsy{B}, T)$,  where $\bsy{B}$ is a partition of $[k]$ and 
 $T$ is a partial tableau of shape $\mu$ such that content$(T) \subseteq \max(\bsy{B})$; here $\max(\bsy{B})$ is the set consisting of the maximum value in each block of $\bsy{B}$.  In this paper, we modify the definition of $\psi$ so that the partial tableau in the image is always standard. 
 To precisely define the map $\psi$ we begin by stating needed definitions and set our notation.

Let $S$ be a finite set. 
A (set) partition of  $S$ is a collection $\bsy{B}=\{B_1, B_2, \dots, B_t\}$ of pairwise disjoint non-empty subsets of $S$ such that $B_1 \cup B_2 \cup \cdots \cup B_k=S$. 
Each $B_i$ is called a block of $\bsy{B}$. Let $\Pi(S)$ be the set of all partitions of $S$. 
For $ 0\leq j \leq k$, 
let  
\begin{eqnarray*} 
\Pi(k,j)&=&\{ \bsy{B}: \ \bsy{B} \in \Pi([k]) \text{ and $\bsy{B}$ has exactly $j$ blocks} \}, \\ 
\Pi^*(k,j)&=&\{ \bsy{B}^*=(\bsy{B}, A): \bsy{B}  \in \Pi([k]) \text{ and  $A$ is a subset of  $j$ blocks of $\bsy{B}$} \}. 
\end{eqnarray*} By definition, if $\bsy{B}^* =(\bsy{B}, A) \in \Pi^*(k,j)$, then $\bsy{B}$ has at least $j$ blocks. 
We say that each block in $A$ is marked and call $\bsy{B}^* \in \Pi^*(k,j)$ \emph{a partition of $[k]$ with $j$ marked blocks}. See Example \ref{ex:1} for an  illustration. 
\vanish{ 
Denote by $B(k,j)$ the cardinality of $\Pi^*(k, j)$. 

Then clearly for $k \geq 1$, 
\begin{equation} \label{B_kj}
  B(k,j) = \sum_{r} \binom{r}{j}  \ S(k,r),    
 \end{equation}
 where $S(k,r)$ is the Stirling number of the second kind, which counts the number of set partitions of $[k]$ into $j$ blocks.
 Set $B(0,0)=1$. 
  The array $B(k,j)$ is given by \href{https://oeis.org/A049020}{A049020} in OEIS.  
} 

Let $\mu \vdash j$.  
The following is the bijection $\psi$ from $\mathcal{SVT}_k(\mu)$ to $\Pi^*(k,j) \times \mathcal{SYT}(\mu)$. 

\medskip 
\noindent 
\underline{The Bijection $\psi$}:   \ \ 
Given a simplified vacillating tableau  ${P} = (\emptyset=\lambda^{(0)}, \lambda^{(\frac{1}{2})}, \lambda^{(1)},\dots, \lambda^{(k)})$, we will recursively define a sequence $(E_0, T_0), (E_{\frac{1}{2}}, T_{\frac{1}{2}}), 
\dots, (E_{k}, T_{k})$,  where for any index $i$, $E_i$  is a set of ordered pairs of integers in $[k]$ (which are viewed as ``edges"), and $T_i$ is a partial tableau of shape $\lambda^{(i)}$.  Let $E_0$  be the empty set and  $T_0$ be the empty tableau. 
For each integer $j=1,2, \dots, k$, assume $(E_{j-1}, T_{j-1})$ is known. 

\begin{enumerate}[(a)] 

   \item If $\lambda^{(j-\frac{1}{2})} = \lambda^{(j-1)}$, then $(E_{j-\frac{1}{2}}, T_{j -\frac{1}{2}}) = (E_{j-1}, T_{j-1})$.  
   
   \item  If $\lambda^{(j-\frac{1}{2})} \subsetneq  \lambda^{(j-1)}$, 
    let $T_{j-\frac{1}{2}}$ be the unique partial tableau with the property that there exists an integer $m$ such that $T_{j-1} = \left(m \xrightarrow{RSK} T_{j-\frac{1}{2}}\right)$.
    Note that $m$ must be less
than $j$. Let $E_{j-\frac{1}{2}}$ be obtained from $E_{j-1}$ by adding the ordered pair $(m, j)$.

  \item If $\lambda^{(j)} =  \lambda^{(j-\frac{1}{2})}$, then 
  $(E_{j}, T_{j})= (E_{j-\frac{1}{2}}, T_{j -\frac{1}{2}})$.  
  \item  If 
   $\lambda^{(j)} \supsetneq  \lambda^{(j-\frac{1}{2})}$, let $E_j = E_{j-\frac{1}{2}}$ 
and $T_j$ be obtained from $T_{j-\frac{1}{2}}$ by adding the entry $j$ in the box  
$ \lambda^{(j)}  / \lambda^{(j-\frac{1}{2})}$. 

\end{enumerate} 
It is clear from the above construction that $E_0 \subseteq E_{\frac{1}{2}} \subseteq \cdots \subseteq  E_{k}$. 
Let $G=(V, E_k)$ be a graph with vertex set $V=[k]$ and edge set $E_k$, and let $\bsy{B}$ be the set partition of $[k]$ whose blocks are vertices of connected components of $G$. 

Let $\mu$ be the shape of $T_k$. Note that if an integer $i$ appears in $T_k$, then $E_k$ cannot contain any ordered pair $(i, j)$ with $i < j$. It follows that $i$ is the maximal element in the block containing it. Hence the content of $T_k$ is a subset of $\max(\bsy{B})$. We get a set partition $\bsy{B}^*$ with marked blocks by putting a mark on each block $X$  if $\max(X)$ is in $T_k$, and then replacing the integers in $T_k$ with integers $1, 2, \dots, j=|\mu|$, following numerical order. 
This results in a SYT $T$ of shape $\mu$. 

Finally we define $\psi({P}) = (\bsy{B}^*, T)$.

\begin{example} \label{ex:1} 
As an example of the map $\psi$, let $k = 7$ and the simplified vacillating tableau be 
\[
\left(\ytableausetup{boxsize=7px} \emptyset\, ,\, \emptyset\, ,\, \ydiagram{1}\, ,\, \ydiagram{1}\, ,\, \ydiagram{1,1}\, ,\, \ydiagram{1,1}\, ,\, \ydiagram{2,1}\, ,\, \ydiagram{1,1}\, ,\, \ydiagram{2,1}\, ,\, \ydiagram{2}\, ,\, \ydiagram{2,1}\, ,\, \ydiagram{2,1}\, ,\, \ydiagram{2,1,1}\, ,\, \ydiagram{2,1}\, ,\, \ydiagram{2,1} \right).
\]

 \ytableausetup{smalltableaux}

Then the sequence of $T_j$ and the corresponding new edge added to $E_j$ at each step are given in the table below. 
\begin{center} 
\begin{tabular}{c|ccccccccccccccc}
   \\
  $j$    &  $0$ & $\frac{1}{2}$  & $1$  & $1\frac{1}{2}$  & $2$ & $2\frac12$ & $3$ & $3\frac12$ & $4$ &  $4\frac12$ & $5$ & $5\frac12$ & $6$ & $6\frac12$ & $7$ \\ \\ 
  \hline 
  \rule{0pt}{1.2\normalbaselineskip} $T_j$  &  $\emptyset$ & $ \emptyset $ 
       & \ytableaushort{1}          & \ytableaushort{1}   
       & \ytableaushort{1,2}    &  \ytableaushort{1,2}   
       &  \ytableaushort{13,2}     &  \ytableaushort{1,2}  
       & \ytableaushort{14,2}     & \ytableaushort{24}  
       &  \ytableaushort{24,5}    &   \ytableaushort{24,5} 
       &   \ytableaushort{24,5,6}  &  \ytableaushort{25,6} 
       &   \ytableaushort{25,6} \\ 
   new edge & & & & & & & & (3,4) & & (1,5) & & & & (4,7) & \\     
\end{tabular}
\end{center} 

It follows  that $E_k=\{ (3,4), (1, 5), (4, 7)\}$, and hence 
\[
\bsy{B}^* =  \{1,5^*\  | \ 2^*\  | \  3,4,7\  | \  6^* \}, \qquad  T=\ytableaushort{12,3} \,. 
\]

\end{example}
 
The map $\psi$ is bijective. This map and its restrictions are our main tool to study  simplified and limiting \vtx. To better understand this map, we describe its inverse. First, given a set partition  of $[k]$,  we represent it by an arc-diagram on the vertex set $[k]$ whose edge set consists of arcs connecting consecutive elements of each block in numerical order. Such a diagram is called a \emph{standard diagram} of the set partition. The figure below gives the standard diagram of the set partition $\{1,5\ | \ 2\  | \ 3,4 ,7 \ |\  6\}$.

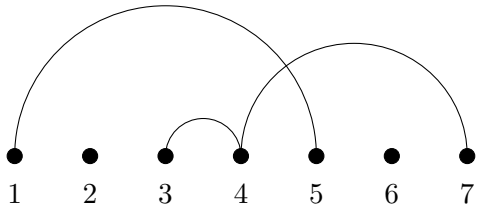
\begin{center}
    \begin{tikzpicture}
        \foreach \r in {1, ..., 7}{
      \filldraw (\r,0) circle (0.1); 
      }
      \foreach \r in {1, ..., 7}{ 
      \node at (\r, -0.5) {$\r$};
      }

      \draw (5,0) arc (0:180:2);  
      \draw (4,0) arc (0:180:0.5);
      \draw (7,0) arc (0:180:1.5); 
      
    \end{tikzpicture}
\end{center}

\noindent 
\underline{The inverse of the map $\psi$}. \ Given $(\bsy{B}^*, T)
\in \Pi^*(k,j) \times \mathcal{SYT}(\mu)$, where $\mu \vdash j$, first we recover the partial tableau $T_k$: Let $a_1 < a_2< \cdots < a_j$ be the maximal elements of the marked blocks in $\bsy{B}^*$. 
Then $T_k$ is obtained from $T$ by replacing entry $i$  with $a_i$ for all $i$. Let 
$G$ be the standard diagram of the set partition of $\bsy{B}^*$.  
We work our way backwards from $T_k$, reconstructing the preceding tableaux and hence the sequence of the shapes, which is the \vt \ $P$. For $j=k, k-1, \dots ,1$, if we know the tableau $T_j$, we can get the tableaux $T_{j-\frac{1}{2}}$ and $T_{j-1}$ by the following rules. 
\begin{enumerate}[(a$^{\prime}$)] 
    \item $T_{j-\frac{1}{2}} =T_j$ if $j$ is not an entry of $T_j$. Otherwise $T_{j-\frac{1}{2}}$ is obtained from $T_j$ be removing the box containing $j$. 
    \item $T_{j-1} =T_{j-\frac{1}{2}}$ if $G$ does not have an edge of the form $(i,j)$ with $i < j$. Otherwise, 
     \[
     T_{j-1}  = \left( i \xrightarrow{RSK} T_{j-\frac{1}{2}} \right). 
     \]
 \end{enumerate}

\begin{example}
    Let $k=3$ and $\lambda=(1, 1)$. There are six simplified \vtx \ 
    of shape $\lambda$ and length $6$. The corresponding pairs
    of $(E_k,T_k)$ and hence $(\bsy{B}^*, T)$ are given in Table \ref{tab:ex}. 
\renewcommand{\arraystretch}{1.5} 
\begin{table}[htb]
\centering
\ytableausetup{smalltableaux}
    \begin{tabular}{|c|c|c|} \hline 
   \vt        & $(E_k, T_k)$  &  $ (\bsy{B}^*, T) $  \\ \hline 
   $ (\emptyset\, , \emptyset\,, \ydiagram{1}\,, \emptyset\,, \ydiagram{1}\,, \ydiagram{1}\,, \ydiagram{1,1} \  )$  & 
        $ (\{(1,2) \}, \   \begin{ytableau}
                           2 \\
                           3   
                         \end{ytableau} \ ) $ 
  &  $ (\{1,2^* \ | \ 3^*\},  \ 
   \begin{ytableau}
     1 \\ 
     2   
   \end{ytableau}  \ ) $  \\ [7pt]\hline 
 $( \emptyset\,, \emptyset\,,  \ydiagram{1}\,, \ydiagram{1}\,, \ydiagram{1,1}\,, \ydiagram{1}\,, \ydiagram{1,1} \ ) $ 
   & $ ( \{(1,3)\}, \   \begin{ytableau}
                          2 \\
                          3   \end{ytableau} \ ) $ 
   & $ (\{1,3^* \ | \ 2^*\},  \ 
   \begin{ytableau}
     1 \\ 
     2   
   \end{ytableau}  \  ) $  \\ [7pt]\hline 

   $( \emptyset\,, \emptyset\,,  \ydiagram{1}\,, \ydiagram{1}\,, \ydiagram{2}\,, \ydiagram{1}\,,\ydiagram{1,1} \ )$ 
   & $ ( \{(2,3)\}, \   \begin{ytableau}
                          1 \\
                          3   \end{ytableau} \ ) $ 
   & $ (\{1^* \ | \ 2, 3^*\},  \ 
   \begin{ytableau}
     1 \\ 
     2   
   \end{ytableau}  \  ) $  \\ [7pt]\hline 
 $( \emptyset\,, \emptyset\,,  \ydiagram{1}\,, \ydiagram{1}\,, \ydiagram{1,1}\,, \ydiagram{1,1}\,, \ydiagram{1,1}\ ) $ 
   & $ ( \emptyset, \   \begin{ytableau}
                          1 \\
                          2   \end{ytableau} \ ) $ 
   & $ (\{1^* \ | \ 2^*\ | \ 3 \},  \ 
   \begin{ytableau}
     1 \\ 
     2   
   \end{ytableau}  \  ) $  \\ [7pt]\hline 

   $( \emptyset\,,\emptyset\,, \emptyset\,, \emptyset\,,   \ydiagram{1}\,, \ydiagram{1}\,, \ydiagram{1,1} \ )$ 
   & $ ( \emptyset, \   \begin{ytableau}
                          2 \\
                          3   \end{ytableau} \ ) $ 
   & $ (\{1 \ | \ 2^*\ | \ 3^* \},  \ 
   \begin{ytableau}
     1 \\ 
     2   
   \end{ytableau}  \  ) $  \\ [7pt]\hline 
  $( \emptyset\,,\emptyset\,, \ydiagram{1}\,, \ydiagram{1}\,,  \ydiagram{1}\,, \ydiagram{1}\,, \ydiagram{1,1} \ ) $ 
   & $ ( \emptyset, \   \begin{ytableau}
                          1 \\
                          3   \end{ytableau} \ ) $ 
   & $ (\{1^* \ | \ 2\ | \ 3^* \},  \ 
   \begin{ytableau}
     1 \\ 
     2   
   \end{ytableau}  \  ) $  \\ [7pt]\hline 
    \end{tabular} 
    \caption{The simplified \vtx\ of shape $\lambda=(1,1)$ and length $6$.}\label{tab:ex}
\end{table} 
    
    \end{example}



\section{Vacillating tableaux  of even length}\label{sec:svt1}

This section pertains to $g_k(\mu)$, the number of simplified \vtx \ of shape $\mu$ and length $2k$. We will remark on instances where a result is applicable to $m_{n,k}^\lambda$ for arbitrary  integers $n\geq 1$ and $k\geq 0$.

\subsection{Formula for \texorpdfstring{$g_k(\mu)$}{}} \mbox{} 

The following theorem is an immediate consequence of the bijection $\psi$. 
\begin{theorem}   \label{g_k-value}
 For all integers $k \geq 0$, 
    \[ g_k(\mu) = B(k,|\mu|) f^\mu,
    \]
     where  $B(k,j)$ is the number of set partitions of $[k]$ with $j$    marked blocks. 
  \end{theorem} 
  It is easy to see that for $0 \leq j \leq k$, 
\begin{equation} \label{B_kj}
  B(k,j) = \sum_{r} \binom{r}{j}  \ S(k,r),    
 \end{equation}
 where $S(k,r)$ is the Stirling number of the second kind that counts the number of set partitions of $[k]$ with exactly $r$ blocks.
 By convention, we let $S(0,0)=1$ and $S(0,r)=0$ for $r \geq 1$. 
The array $B(k,j)$ is given by \href{https://oeis.org/A049020}{A049020} in OEIS\cite{OEIS}.  
 The exponential generating function of $B(k,j)$ for a fixed $j$ is  
 \begin{eqnarray} \label{egf:b}
\sum_{k \geq 0} B(k,j) \frac{x^k}{k!}  = \frac{1}{j!}(e^x -1)^j \exp(e^x -1).
\end{eqnarray}   

 Theorem \ref{g_k-value} and  Equation \eqref{egf:b} are given in \cite[Thm 2.4]{CDDSY07} and Formula \eqref{B_kj} is given  in \cite[Eq.(5.11)]{BH17}. 
When $\mu=\emptyset$, Theorem \ref{g_k-value} implies $g_k(\emptyset)=B(k)$, the $k$-th Bell number that counts the number of set partitions of $[k]$.  Further refinement between simplified vacillating tableaux of shape $\emptyset$ and set partitions was studied in \cite{CDDSY07}. 

 \textsc{Remark}.    The special case $\mu=\emptyset$ of Theorem \ref{g_k-value}, written in terms of
 $m_{n,k}^\lambda$, gives $m_{n,k}^{(n)} =B(k)$ for $n \geq 2k$. 
  Martin and Rollet \cite{MR98} proved a more general identity 
\begin{equation} \label{vt-case1}
\sum_{j=1}^n S(k, j) = m_{n,k}^{(n)}, \qquad \text{ for }  n,k  \geq 1. 
\end{equation}
Equation \eqref{vt-case1} with  even $k$  is also proved by Benkart and Halverson in \cite{BH17}. 
Using growth diagrams, 
Krattenthaler \cite{Kratten23}   gave a  simple proof of \eqref{vt-case1} as well as the identity 
$S(k,n)+S(k,n-1)= m_{n,k}^{(1^n)}$ for $n \geq 1$. 
Note that the latter identity is trivial when $n \geq 2k$. 

Using the same construction as $\psi$ but starting with a Young tableau of shape $(n)$ filled with $n$ zeros, Benkart, Halverson and Harman in \cite{BHH17} gave a bijection from the set of $n$-\vtx \ in $\mathcal{VT}_{n,k}(\lambda)$ to pairs 
$(\bsy{B}, \hat{T})$, where $\bsy{B} \in \Pi(k,j)$ is a partition of $[k]$ into $j$ blocks, and $\hat{T}$ is a SSYT of shape $\lambda$ with  content 
$\{0^{n-j}\} \cup \max(\bsy{B})$.  We will use
 this result in Subsections \ref{ssec:schur2} and \ref{ssec:schur2''} to get the expansion of $\sum_{\lambda} m_{n,k}^\lambda s_\lambda(\bsy{x})$
in terms of complete symmetric functions. 
See Theorems \ref{thm:gk-schur2} and \ref{thm:gk'-schur2}.

\subsection{The sum \texorpdfstring{$\sum_\mu g_{k_1}(\mu) g_{k_2}(\mu)$ }{}}\mbox{} 

Note that there is a symmetry between the even steps and the odd steps in the definition of simplified vacillating tableaux. 
Thus any walk on  Young's lattice, as described in Definition \ref{def:svt}, from the empty shape to itself in $2(k_1+k_2 )$ steps can be viewed as a walk from $\emptyset$  to some shape $\mu$  in $2k_1$  steps, then followed by the reverse of a walk from $\emptyset$ to $\mu$ in $2k_2$ steps.  It follows that
for all integers $k_1, k_2 \geq 0$, 
\begin{equation} \label{g^2}
    \sum_\mu g_{k_1}(\mu) g_{k_2}(\mu) = g_{k_1+k_2}(\emptyset) = B(k_1+k_2). 
\end{equation}
For the special case where $k_1=k_2$, Identity \eqref{g^2} is proved in \cite{HL05} in the form $\sum_{\lambda \in \Lambda_n^k} (m_{n,k}^\lambda)^2=B(2k)$ for $n\geq 2k$. 

\subsection{The sum  \texorpdfstring{$g_k=\sum_\mu g_k(\mu)$}{}} \mbox{} 

Let $g_k=\sum_\mu g_k(\mu)$ be the number of simplified vacillating tableaux of length $2k$.  
The initial  terms of the sequence $(g_k)_{k=0}^\infty$ are $1, 2, 7, 31, 164, 999, \dots$. This sequence is cataloged as \href{https://oeis.org/A002872}{A002872} in OEIS. 

Let $[-k]$ be the set of integers $\{-k, \ldots, -2, -1\}$. 
A set partition of $[-k] \cup [k]$ is \emph{symmetric} if $-B$ is a block whenever $B$ is  a block of the set partition.  The following theorem was first proved in \cite[Eq.(5.5)]{HL05}, using the symmetry of the bijection 
$\psi$ when one restricts $\psi$ to a  simplified \vt  \ of shape $\emptyset$ and its reverse. 
Here we give a purely 
combinatorial proof  using the pairs $(\bsy{B}^*, T)$  in  $\Pi^*(k,j) \times \mathcal{SYT}(\mu)$, where $|\mu|=j$.

\begin{theorem} \label{thm:vk}
Let $k \geq 1$ be an integer. Then $g_k$ counts the number of symmetric partitions of $[-k] \cup [k]$. 
\end{theorem} 
\begin{proof} \  
 Each simplified \vt \ of length $2k$ is uniquely represented by a pair $(\bsy{B}^*, T) \in 
 \Pi^*(k,j) \times \mathcal{SYT}(\mu)$  with $|\mu|=j$. 
 Draw the standard diagram of the partition $\bsy{B}$ on vertices 
 $1, 2, \dots, k$, followed by the reverse of this diagram on vertices with labels $-k, \dots, -2, -1$. 
 If the marked blocks of $\bsy{B}^*$ are $X_1, \dots, X_j$, (say, ordered by the values of their maximal elements), then the corresponding marked blocks on $[-k]$  are $-X_1, \dots, -X_j$. 
 For the SYT $T$, applying the inverse of RSK to the pair $(T, T)$ we obtain an involution $\sigma \in \mathfrak{S}_j$. 
 Now connect the maximal element of block $X_i$ with the minimal element of block $-X_{\sigma(i)}$, for each $i=1, \dots, j$. 
Since $\sigma$ is an involution, it means that the maximal element of block $X_{\sigma(i)}$ is connected to the minimal element 
 of $-X_i$, and hence  we get a symmetric diagram that defines a symmetric set partition of $[-k] \cup [k]$. 
Note that  the marked blocks of $\bsy{B}^*$ correspond to the blocks of the symmetric set partition with both positive and negative integers. 
 Each step of the above construction  can be easily reversed. Hence our construction is a bijection.
 \vanish{ 
 In the proof of \eqref{g^2}, let one walk from the empty shape to some  shape $\mu$ in $2k$ steps, then 
 reverse the same walk back to the the empty shape.  The number of such walks is counted by $g_k$.  
 Via the bijection $\psi$  each such a walk (from $\emptyset$ to $\emptyset$) corresponds to a set partition on $[2k]$ whose standard diagrams are symmetric with respect to the middle line.  Then by \cite[Theorem 5.5]{HL05} 
 using growth diagram that $\psi$ restricted to is symmetric, in the sense that symmetric 
 
 Re-labeling the elements of $[2k]$ by $\{-k, \cdots, -1, 1, \dots, k\}$  from left to right, we  get a symmetric set partition of $[-k] \cup [k]$.  
 } 
\end{proof}

\begin{example}
    As an example, let $B^*=\{ 1^*\ |\ 2,4\ | \ 5^*\ |\ 3, 6^*\}$ and $T= \begin{ytableau}
1 & 3 \\
2
\end{ytableau}$.   The marked blocks are ordered by their maximal elements, therefore $X_1=\{1\}$, $X_2=\{ 5\}$, and $X_3=\{3,6\}$. 
The involution determined by $T$ is  $\sigma=213$ (in one-line notation). Hence we merge the blocks $X_1$ with $-X_2$, $X_2$ with $-X_1$, 
and $X_3$ with $-X_3$, as shown in the next figure. The resulting symmetric set partition is $\{ 1, -5\  |\  5, -1\ | \ 2,4\  |\  -2, -4\  | \ 3,6, -3, -6\}$
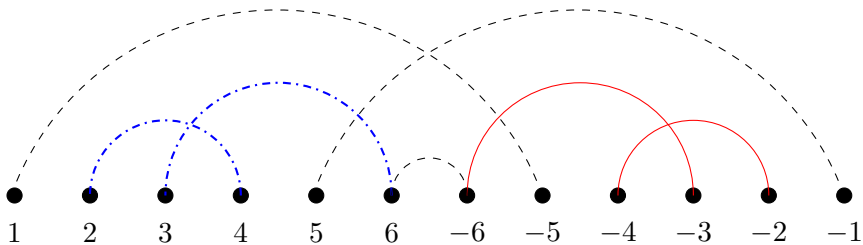
\begin{center}
    \begin{tikzpicture}
     \foreach \r in {1, ..., 12}{
      \filldraw (\r,0) circle (0.1); 
      }
      \foreach \r in {1, ..., 6}{ 
      \node at (\r, -0.5) {$\r$};
      \node at (13-\r, -0.5) {$-\r$}; 
      }
      
       \draw[blue, thick, dash dot] (4,0) arc (0:180:1);  
      \draw[blue, thick, dash dot ] (6,0) arc (0:180:1.5); 
  
      \draw[red]  (10,0) arc (0:180: 1.5); 
   \draw[red] (11,0) arc (0:180:1); 
   
      \draw[dashed] (8,0) arc (20:160: 3.75); 
      \draw[dashed] (12,0) arc (20:160:3.75); 
      \draw[dashed] (7,0) arc (0:180:0.5); 
    \end{tikzpicture}
\end{center}    
\end{example}

\subsection{The sum  \texorpdfstring{$\sum_\mu g_k(\mu) f^\mu$}{}} \mbox{}

Similar to Identity \eqref{n-vt}, let us consider  the sum $u_k =\sum_\mu g_k(\mu) f^\mu$. 
The initial values of $(u_k)_{k=0}^\infty$ are $1, 2, 7, 33, 198, \ldots$, which coincide with the initial terms of \href{https://oeis.org/A059099}{A059099} in OEIS that was studied by Nkonkobe and Murali  \cite{NM15}
as one case of ``restricted barred preferential arrangements.'' 
We will present a simple combinatorial proof that $(u_k)_{k=0}^\infty$ is indeed the sequence \href{https://oeis.org/A059099}{A059099}.   
  
Recall that 
an \emph{ordered partition}, or the \emph{preferential arrangement}  of a set $S$,  is a partition of $S$ into disjoint non-empty blocks, together with a linear order on the blocks. 
For $n \in \mathbb{N}$, ordered partitions of $[n]$ are counted by  \emph{Fubini numbers}, which  
 have exponential generating function  $1/(2-e^x)$ and 
appear as \href{https://oeis.org/A000670}{A000670} in OEIS. 

We say that a set partition of $S$ is \emph{partly ordered}  if it has a (possibly empty) subset of marked blocks that are linearly ordered.
For example, the following set partition of $\{1, 2, \dots, 7\}$  is  partly ordered:
\[
\{  \left( 2,7^* \ | \   1,3^*\ \right) \ |   \ 4,6 \ | \  5 \}.  
\]
The first two blocks with marks are linearly ordered
(hence listed inside the parentheses); the last two blocks are not marked and, hence, are  unordered. 
If we switch the first two blocks, we obtain a different partly ordered set partition, while switching the last two blocks gives the same partly ordered set partition. That is,
\[
\{  \left( 1,3^* \ | \   2,7^*\ \right) \ |   \ 4,6 \ | \  5 \} \neq \{  \left( 2,7^* \ | \   1,3^*\ \right) \ |   \ 4,6 \ | \  5 \} = \{  \left( 2,7^* \ | \   1,3^*\ \right) \ | \ 5 \ |  \ 4,6  \} .  
\]

\begin{theorem} \label{thm:pk}
    Let $k\geq 1$ be an integer. The number $u_k$ counts the set of partly ordered set partitions of $[k]$. 
\end{theorem}
\begin{proof}
    Again we use the bijection $\psi$ to represent a simplified vacillating tableau  of length $2k$ by  a pair
    $(\bsy{B}^*, T)$, where $\bsy{B}^*$ is a set partition of $[k]$ with $j$ marked blocks,  and $T$ is a SYT of some shape $\mu$ with $|\mu|=j$.  Now $u_k$ counts the pairs in $\mathcal{SVT}_k(\mu) \times \mathcal{SYT}(\mu)$, 
    each of which corresponds to a triple $(\bsy{B}^*, T, S)$ where both $T$ and $S$ are SYTs of shape $\mu$. Via the inverse of the RSK algorithm, $(T,S)$ uniquely determines a permutation $\sigma \in \mathfrak{S}_j$, which gives the linear order 
    on the marked blocks of $\bsy{B}^*$. 
\end{proof}
From Theorem \ref{g_k-value} we have 
$u_k=\sum_j j! B(k,j)$, where $B(k,j)$ is given by Formula \eqref{B_kj}. 
By definition of partly ordered set partition,  the exponential generating function of $(u_k)_{k=0}^\infty$ can be expressed as the product of the exponential generating functions of Bell numbers and Fubini numbers. Explicitly, 
\begin{equation} \label{egf_pk}
\sum_{k \geq 0} u_k \frac{x^k}{k!}  =  \frac{\exp(e^x-1) }{2-e^x}. 
\end{equation}

\subsection{Product of \texorpdfstring{$g_k(\mu)$}{} with Schur functions \texorpdfstring{$\sum_\mu g_k(\mu) s_\mu(\bsy{x})$}{}}  \mbox{}

Next we connect vacillating tableaux to symmetric functions. Let $\mu$ be a partition. For a SSYT $T$ of shape $\mu$ and content$(T)\subseteq \mathbb{Z}^+$, let $\alpha_i(T)$ be 
the number of times the integer $i$ appears in $T$, and write
\[
\bsy{x}^{T} = x_1^{\alpha_1(T)} x_2^{\alpha_2(T)} \cdots. 
\]
The Schur function $s_\mu$ in the variables $\bsy{x}=(x_1, x_2, \dots)$ is the formal power series 
\[
s_\mu (\bsy{x}) = \sum_T \bsy{x}^T,
\]
where $T$ ranges over all SSYT of shape $\mu$ and content$(T)\subseteq \mathbb{Z}^+$.  By convention, set 
$s_\emptyset =1$. 

\begin{theorem} \label{Thm:gk-schur1} 
Let  $k \geq 0$ be an integer.  
We have the  identity 
    \begin{equation} \label{gk-schur1}
         \sum_\mu g_k(\mu) s_\mu(\bsy{x} ) = \sum_{j=0}^k B(k, j) h_1(\bsy{x})^j,
    \end{equation}
    where $B(k,j)$ is  the number of set partitions of $[k]$ with $j$  marked blocks and $h_1(\bsy{x})=\sum_i x_i$.  
\end{theorem}
\begin{proof}
    The left side of \eqref{gk-schur1} is a summation of $\bsy{x}^S$ over the set of triples $(\bsy{B}^*, T, S)$ where $\bsy{B}^\ast$ is a set partition of $[k]$ with  $j$ marked blocks, $T$ is a SYT of some shape $\mu$ of size $j$, and $S$ is a  SSYT of shape $\mu$ and content$(S)\subseteq \mathbb{Z}^+$.   The pair $(T, S)$, under the inverse of  RSK algorithm, corresponds uniquely to a two-line array of the form \eqref{2-line-array}, where the top line 
    is $(1, 2, \dots, j)$, and the second line is a sequence $\bsy{i}=(i_1, i_2, \dots, i_j)$ of positive integers. 
     Note that  
    \[
    \bsy{x}^S = x_{i_1} x_{i_2} \cdots x_{i_j}. 
    \]
    Hence 
    \[
    \sum_{(\bsy{B}^*, T, S)} \bsy{x}^S = \sum_j B(k, j) 
    \sum_{\bsy{i} \in (\mathbb{Z}^+)^j} x_{i_1} x_{i_2} \cdots x_{i_j} 
    = \sum_j B(k,j) \left( \sum_i x_i\right)^j.\qedhere
    \]    
\end{proof}

Let $m$ be a positive integer and $x_i=q^{i-1}$ for 
$i=1, 2, \dots, m$ and $x_j=0$ for $j >m$. 
Then  Identity \eqref{gk-schur1} becomes 
\begin{equation} \label{gk-Schur-q} 
\sum_{\mu}  g_{k}(\mu)  s_\mu(1, q, \dots, q^{m-1}) 
=\sum_{j=0}^k B(k,j)  [m]_q^j, 
\end{equation}
where $[m]_q=1+q+\cdots +q^{m-1}$ is the $q$-integer. 
The polynomial $s_\mu(1, q, \dots, q^{m-1})$ is called the principal specialization of $s_\mu$, 
which can be computed by Stanley's hook-content formula \cite[Theorem 7.21.2]{EC2}. For an integer partition $\mu=(\mu_1, \dots, \mu_t)$, let $\mu' = (\mu'_1,..., \mu'_{\mu_1})$ be the conjugate of $\mu$, where $\mu'_j$ is the number of boxes in the $j$-th column of the Young's diagram of $\mu$.   
For a box $u=(i,j)$ in the $i$-th row and $j$-th column of the Young diagram of 
$\mu$, define the \emph{hook length} $h(u)$
of $\mu$ at $u$  by 
\[
h(u)=\mu_i+ \mu'_j -i-j+1. 
\]
In addition, define 
\[
b(\mu) = \sum_{i \geq 1} (i-1)\mu_i = \sum_{i \geq 1}  \binom{\mu'_i}{2}. 
\]
Then  the hook-content formula is  
 \begin{equation} \label{q-hook-content}
     s_\mu(1, q, \dots, q^{m-1}) = q^{b(\mu)}
     \prod_{ u = (i,j) \in \mu} \frac{[m+j-i]_q}{[h(u)]_q}.
 \end{equation}

\subsection{Product  of \texorpdfstring{$m_{n,k}^\lambda$}{} with Schur functions \texorpdfstring{$\sum_\lambda  m_{n,k}^\lambda s_\lambda(\bsy{x} ) $}{}}  \mbox{} 
\label{ssec:schur2}

To compare with the result of the preceding subsection, we also consider the product of Schur functions with the number of $n$-\vtx\  for general $n\geq 1$ and $k\geq 0$. 
Note that this summation is over all integer partitions $\lambda \in \Lambda_n^k$, while for $g_k(\mu)$, $\mu$ ranges over all integer partitions of size at most $k$. 

For a nonnegative integer $n$, the complete symmetric functions $h_\lambda$ are defined by the formulas
\begin{eqnarray*}
h_0&=&1, \\
    h_n& = &\sum_{i_1 \leq  \cdots \leq i_n} x_{i_1} \cdots x_{i_n}, \qquad  (n \geq 1) \\ 
    h_\lambda & = & h_{\lambda_1} h_{\lambda_2} \cdots  \qquad \text{ if } \lambda=(\lambda_1, \lambda_2, \dots ).  
\end{eqnarray*}
The sum of products of $m_{n,k}^\lambda$ and Schur functions has a nice expansion in complete symmetric functions.

\begin{theorem} \label{thm:gk-schur2} For integers 
$k \geq 0$ and $n \geq 1$,
    \begin{eqnarray*} \label{gk-schur2}
        \sum_{\lambda \in \Lambda_n^k}  m_{n,k}^\lambda  s_\lambda(\bsy{x} )
        &=& \sum_{j=0}^{\min(n,k)} S(k,j) h_{(n-j, 1^j)}(\bsy{x}),        \end{eqnarray*}
        where if $j=n$, the integer partition $(n-j, 1^j)$ is just $(1^n)$. 
\end{theorem}
\begin{proof}
   We use a bijection of Benkart, Halverson, and Harman \cite{BHH17} that maps 
    the set of $n$-vacillating tableaux of shape $\lambda$ and length $2k$  to the set of pairs $(\bsy{B}, \hat{T})$ where $\bsy{B} \in \Pi(k,j)$ is a partition of $[k]$ into $j$ blocks   for some $j \leq n$, and $\hat{T}$ is a SSYT of shape $\lambda$ filled with $n-j$ zeros and $j$ distinct positive integers, which are  the maximal elements of the blocks of $\bsy{B}$. 

   We extend the above bijection to an injection from $\mathcal{VT}_{n,k}(\lambda)$ into a subset of $\Pi([k]) \times \mathcal{SYT}(\lambda)$ and then analyze the images. 
   First, we replace the SSYT $\hat{T}$ with a SYT $T$ of shape $\lambda$ as follows: 
   Replace the $n -j$ zeros, all of which must be in the first row of $\hat{T}$, by $1,2, \ldots, n-j$ from left to right. 
   Then replace the remaining $j$ entries of $\hat{T}$ by $n-j+1, \dots, n$ according to their numerical order. Conversely,  given $\bsy{B}=\{B_1, \dots, B_j\} $ and $T$,  we can recover $\hat{T}$ whose content is $\{0^{n-j}\} \cup \{ \max(B_i): i=1, 2, \dots, j\}$ simply by replacing each of  $1, \dots, n-j$ in $T$ with $0$ and replacing $n-j+i$ in $T$ with the $i$-th smallest element in $\max(\bsy{B})$.
   
    Now we can represent an $n$-\vt\ by a pair $(\bsy{B}, T)$ where $\bsy{B}$ is a partition of $[k]$ with $j$ blocks, and $T$ is a SYT of some shape $\lambda \vdash n$ with the property that $\lambda_1 \geq n-j$ and the entries  $1, \dots, n-j$ are in the first row of $T$. 
    
    The formula $\sum_{\lambda \in \Lambda_n^k}  m_{n,k}^\lambda  s_\lambda(\bsy{x})$ 
   is the formal power series that sums $\bsy{x}^S$ over the set of triples $(\bsy{B}, T, S)$ where $(\bsy{B}, T)$ is as described above,  and $S$ is a SSYT of shape $\lambda$ with positive integer entries. Fix the set partition $\bsy{B} \in \Pi(k,j)$ and consider all the tableaux $T$ and $S$ that can appear with $\bsy{B}$ in some triples. 
   Under the inverse of RSK insertion, taking $S$ as the insertion tableau and $T$ as the recording tableau, the pair $(S, T)$ uniquely corresponds to a  sequence 
   $\bsy{i} =(i_1, \dots, i_n)$ of positive integers, where $\bsy{x}^S= x_{i_1} \cdots x_{i_n}$. 
   Since $T$ is the recording tableau, 
   the condition  that the entries $1, \dots, n-j$ appear in the first row of $T$  holds if and only if    $i_1 \leq \cdots \leq i_{n-j}$. 
   Let 
    \[
   \mathcal{D}_{n, n-j}=\{ (i_1, \dots, i_n) \in (\mathbb{Z}^+)^n: \ i_1 \leq \cdots \leq i_{n-j} \}. 
\]
 Then the formal power series 
 $ \sum_{\lambda \in \Lambda_n^k}  m_{n,k}^\lambda  s_\lambda(\bsy{x} ) $  
 can be expressed as 
  \begin{eqnarray*}
  \sum_j S(k,j) \sum_{\bsy{i}\in \mathcal{D}_{n,n-j}} x_{i_1} \cdots x_{i_n} 
 & = &\sum_j S(k,j) \left( \sum_{i_1 \leq \cdots \leq i_{n-j} } x_{i_1}\cdots x_{i_{n-j}} \right) \left(\sum_{i} x_i \right)^j \\ 
 & =& \sum_{j} S(k,j) h_{(n-j, 1^j)}(\bsy{x} ),
  \end{eqnarray*} 
  where the range of $j$ is $0 \leq j \leq \min(n,k)$.
\end{proof}

\begin{example}
     We illustrate Theorem \ref{thm:gk-schur2} with some examples. 
     First, let $n=2$ and $k \geq 1$.  From \cite{Kratten23} we have that 
 $m_{n,k}^{(2)}=m_{n,k}^{(1,1)}=S(k,1)+S(k,2)=2^{k-1}$.
    Hence Theorem \ref{thm:gk-schur2} becomes 
    $$ 
     \sum_{\lambda \in \Lambda_2^k}  m_{2,k}^\lambda  s_\lambda(\bsy{x} ) = 2^{k-1} (s_{(2)} + s_{(1,1)} ) = 2^{k-1} h_{(1,1)}
     =S(k,1) h_{(1,1)} + S(k,2) h_{(1,1)}.
     $$ 
     As another example, let $n=6$ and $k=3$. Then we have 
     
     \begin{eqnarray*}
       \sum_{\lambda \in \Lambda_6^3}  m_{6,3}^\lambda  s_\lambda(\bsy{x}) &=& 5 s_{(6)} + 10s_{(5,1)} + 6s_{(4,2)} + 6s_{(4,1,1)}+s_{(3,3)} + 2s_{(3,2,1)} + s_{(3,1,1,1)}\\
       &=& h_{(5,1)} + 3h_{(4,1,1)} + h_{(3,1,1,1)}\\ 
       \\
       &=& S(3,1)h_{(5,1)} + S(3,2)h_{(4,1,1)} + S(3,3)h_{(3,1,1,1)}.  
     \end{eqnarray*}
\end{example}

For the principal specialization, note that   
$h_n=s_{(n)}$. Using Stanley's hook-content formula, we have 
$h_n(1, q, \dots, q^{m-1})= \qbinom{m+n-1}{n}$. Therefore, Theorem \ref{thm:gk-schur2}  implies 
\begin{equation} \label{Schur-q2} 
\sum_{\lambda \in \Lambda_n^k}  m_{n,k}^\lambda  s_\lambda(1, q, \dots, q^{m-1}) 
=\sum_{j \leq \min(k,n) } S(k,j) \qbinom{m+n-j-1}{n-j} [m]_q^j. 
\end{equation}

%
%

\section{Vacillating tableaux of odd length}
\label{sec:svt2}

\subsection{Formula for \texorpdfstring{$g_{k+\frac{1}{2}}(\mu)$}{}} \mbox{}  

In the description of the bijection $\psi$  in Subsection \ref{ssec:vt}, comparing the pairs  $(E_i, T_i)$ for the last two steps where $i=k-\frac{1}{2}$ and $i=k$, one observes that $E_{k-\frac{1}{2}}=E_k$ and $T_{k-\frac{1}{2}}$ is either the same as $T_k$ or  missing one corner box containing $k$.  That is, without the last pair $(E_k, T_k)$, we still get the set partition $\bsy{B}$ of $[k]$
with some $j$ blocks marked, and a SYT $T=T_{k-\frac{1}{2}}$ of shape $\mu$ with $|\mu|=j$, with the additional condition that  
the block containing element $k$ is not marked.  

Let $\tilde \Pi^*(k,j) = \{\bsy{B}^*=(\bsy{B}, A): \bsy{B} \in \Pi([k]),\  A  
\text{ is a subset of $j$ blocks of $\bsy{B}$  not containing $k$}\}$. 
Denoted by $\tilde B(k,j)$ the cardinality of $\tilde \Pi^*(k,j)$. 
Then the map $\psi$ induces a bijection 
\[
\mathcal{SVT}_{k-\frac{1}{2}}(\mu) \longleftrightarrow 
\tilde \Pi^*(k,|\mu|) \times \mathcal{SYT}(\mu), 
\]
where  $\mathcal{SVT}_{k-\frac{1}{2}}(\mu) $ is the set of simplified \vtx\ of shape $\mu $ and length $2k-1$. 
It follows that 
 $g_{k-\frac{1}{2}}(\mu) = \tilde B(k, |\mu|) f^\mu$,   or equivalently, 
\begin{equation} \label{tilde-B}
    g_{k+\frac{1}{2}}(\mu) = \tilde B(k+1, |\mu|) f^\mu, 
\end{equation}
where for $0 \leq j \leq k$, 
\begin{equation} \label{B'_kj}
\tilde B(k+1, j) = \sum_r  \binom{r}{j} S(k+1,r+1). 
\end{equation} 
The array $\tilde B(k, j)$ is given in \href{https://oeis.org/A137597}{A137597} in OEIS. 
Formula \eqref{tilde-B} appears in \cite[Eq.(5.12)]{BH17}. 

We remark that when $\mu=\emptyset$, Formula~\eqref{tilde-B} gives that $g_{k+\frac12}(\emptyset)=B(k+1)$, the number of set partitions of $[k+1]$. 

\subsection{The sum  \texorpdfstring{$\sum_\mu g_{k_1+\frac{1}{2}}(\mu) g_{k_2+\frac{1}{2}}(\mu)$}{} } \mbox{} 

\begin{theorem} \label{thm:oddg}
For all integers $k_1,k_2\geq 0$, 
    \begin{eqnarray} \label{odd-gxg}
        \sum_\mu g_{k_1+\frac{1}{2}}(\mu) g_{k_2+\frac{1}{2}}(\mu) =B(k_1+k_2+1). 
    \end{eqnarray}
\end{theorem}
\begin{proof}[First proof of Theorem \ref{thm:oddg}]\ 
Observe that if a simplified \vt \ of shape $\mu$ and length $2k_1+1$ is followed by the reverse of 
another simplified \vt\ of same shape $\mu$ and length $2k_2+1$, 
we obtain exactly a simplified \vt \  of the empty shape $\emptyset$ and length $2k_1+2k_2+2$. Since $g_{k_1+k_2+1}(\emptyset)=B(k_1+k_2+1)$,  Identity \eqref{odd-gxg} follows. 
\end{proof} 

In addition to the above simple argument, we will give another proof that uses the pair $(\bsy{B}^*, T)$ in $\tilde \Pi^*(k,j) \times \mathcal{SYT}(\mu)$. This new argument can be applied to compute  $\sum_\mu g_{k_1}(\mu) g_{k_2+\frac{1}{2}}(\mu)$ and 
$\sum_\mu g_{k_1}(\mu) a_{k_2}(\mu)$, etc. 

\begin{proof}[Second proof of Theorem \ref{thm:oddg}]\ 
    The left  side of \eqref{odd-gxg}  counts the quadruple $(\bsy{B}^*_1, T_1; \bsy{B}^*_2, T_2)$ where for $i=1, 2$, 
    $T_i$ is a SYT of the shape $\mu$ with $|\mu|=j$ for some integer $j$, and 
  $\bsy{B}^*_i$ is a partition of $[k_i+1]$ with $j$ marked blocks such that the block containing $k_i+1$ is not marked. Order the marked blocks of $\bsy{B}^*_i$ by their maximal elements, and let $\sigma \in \mathfrak{S}_j$ be the permutation corresponding to the pair $(T_1, T_2)$ under 
  the RSK algorithm. 

   Draw the standard diagram of the partition $\bsy{B}^*_1$, followed by the reverse of the standard diagram of $\bsy{B}^*_2$.  Then, for each $i \leq j$, 
   add an edge connecting the maximal element of the $i$-th marked block of $\bsy{B}^*_1$ to the maximal element of the 
   $\sigma(i)$-th marked block of $\bsy{B}^*_2$, and then identify the vertex $k_1+1$ of $\bsy{B}^*_1$ with the vertex $k_2+1$ of $\bsy{B}^*_2$.  Now we get a diagram on $k_1+k_2+1$ vertices,  which corresponds to  a set partition of $[k_1+k_2+1]$. The procedure can be easily reversed, and hence it is a bijection.  
\end{proof} 

\begin{example}
     This example illustrates the construction in the second proof. Let  $k_1=7$ and $k_2=6$.
     Take $\bsy{B}_1^*=\{2^*\ | \ 1,3\ | \ 6^*\ |\  5,7^*\ | 4,8 \}$. To distinguish the elements, we put a bar on each element in the partition $\bsy{B}^*_2$. Let 
     $\bsy{B}^*_2= \{ \bar{1}, \bar{4}^*\ | \ \bar{5}^*\ | \ \bar{3}, \bar{6}^*\ | \  \bar{2}, \bar{7}\}$, and 
     \[
     T_1= \begin{ytableau}
1 & 3 \\
2
\end{ytableau}   \quad , \quad
T_2= \begin{ytableau}
1 &2 \\
3
\end{ytableau}\ . 
     \]
The marked blocks are ordered by their maximal elements.  Hence for $\bsy{B}^*_1$, the marked blocks are $X_1=\{2\}$, $X_2=\{6\}$ and $X_3=\{5,7\}$; for $\bsy{B}^*_2$, 
the marked blocks are $Y_1=\{\bar{1}, \bar{4}\}$, $Y_2=\{\bar{5}\}$ and $Y_3=\{ \bar{3}, \bar{6}\}$. 
The permutation determined by $(T_1, T_2)$ is $\sigma=231$ (in one-line notation). 
Hence we merge the blocks $X_1$ with $Y_2$, 
$X_2$ with $Y_3$, and $X_3$ with $Y_1$, and identify 
the element $8$ of $\bsy{B}^*_1$ with the element $\bar{7}$ of $\bsy{B}^*_2$, as shown in the figure below. 
Finally replacing the integer $\bar{i}$ with $k_1+k_2+2-i$,  we get the partition $\{ 1,3\  | \ 2,10\ | \ 4,8, 13\ | \ 5,7,11, 14 \ | \ 6,9, 12 \} $. 

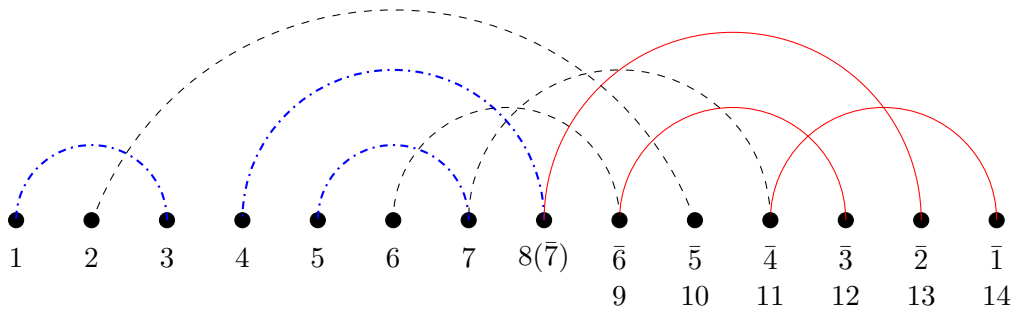
\begin{center}
    \begin{tikzpicture}
     \foreach \r in {1, ..., 14}{
      \filldraw (\r,0) circle (0.1); 
      }
      \foreach \r in {1, ..., 6}{ 
      \node at (\r, -0.5) {$\r$};
      \node at (15-\r, -0.5) {$\bar{\r}$}; 
      }
      \foreach \r in {9, ..., 14}{
       \node at (\r, -1) {$\r$}; 
      }
       \node at (8, -0.5) {$8(\bar{7})$}; 
       \node at (7, -0.5) {$7$}; 
      \draw[color=blue,thick,dash dot] (3,0) arc (0:180:1);  
      \draw[blue, thick, dash dot] (8,0) arc (0:180:2); 
      \draw[blue,thick, dash dot] (7,0) arc (0:180:1); 
      \draw[red]  (14,0) arc (0:180: 1.5); 
      \draw[red] (13,0) arc (0:180:2.5); 
      \draw[red] (12,0) arc (0:180:1.5); 
      \draw[dashed] (10,0) arc (20:160: 4.25); 
      \draw[dashed] (9,0) arc (0:180:1.5); 
      \draw[dashed] (11,0) arc (0:180:2); 
    \end{tikzpicture}
\end{center}     
\end{example}

The special case of \eqref{odd-gxg} with  $k_1=k_2$ was obtained in \cite[Eq.(5.1)]{HL05}. 

With a similar argument as in the second proof of Theorem \ref{thm:oddg}, but without identifying the last elements 
of the partitions  $\bsy{B}_1^*$ and $\bsy{B}^*_2$, we obtain the following result.
\begin{theorem}
  For all integers $k_1, k_2 \geq 0$, 
    the number $\sum_\mu g_{k_1+\frac{1}{2}}(\mu) g_{k_2}(\mu)$ counts the set of partitions 
of $[k_1+k_2+1]$ with the condition that $k_1+1$ is the maximal element in the block containing it.
\end{theorem}

\subsection{The sum  \texorpdfstring{$\sum_\mu g_{k+\frac{1}{2}}(\mu)$}{}} \mbox{}

Let $g_{k+\frac{1}{2}}$ be the sum $\sum_\mu g_{k+\frac{1}{2}}(\mu)$. We will prove in Theorem \ref{thm: number of odd length SVT}
that $g_{k+\frac{1}{2}}$ is the same as the number of symmetric partitions of the set $[-k,k]=\{-k, \dots, -1, 0, 1, \dots, k\}$.
The sequence $(g_{k+\frac{1}{2}})_{k=0}^\infty$ is  \href{https://oeis.org/A080337}{A080337} in OEIS, with the initial values $1, 3, 12, 59, 339\dots$. 

Indeed, combining with the sequence \href{https://oeis.org/A0002872}{A002872} and considering $g_0, g_{\frac{1}{2}}, g_1, g_{1\frac{1}{2}}, g_2, \dots$ , we get the sequence  \href{https://oeis.org/A080107}{A080107} in OEIS, which is described as the number of fixed points of set partitions under the involution  $i \leftrightarrow n+1-i$. It is clearly equivalent to symmetric set partitions.
 
 \begin{theorem}\label{thm: number of odd length SVT}
    For all integers $k \geq 0$, 
    $ g_{k+\frac{1}{2}}$ counts the number of symmetric partitions of the set 
    $[-k,k]=\{-k, \dots, -1, 0, 1, \dots, k\}$. 
\end{theorem}
\begin{proof}
    As in the second proof of Theorem \ref{thm:oddg}, 
    for a pair $(\bsy{B}^*, T)$ counted by $g_{k+\frac{1}{2}}(\mu)$, draw the standard diagram of the partition $\bsy{B}^*$ on vertices $1, 2, \dots, k+1$, followed by the reversed of this diagram on vertices with labels 
    $(k+1)', \dots, 2', 1'$ from left to right.   If the marked blocks of $\bsy{B}^*$ are $X_1, \dots, X_j$,   then the corresponding marked blocks on $[(k+1)']$ are $X_1', \dots, X_j'$. The SYT $T$ corresponds uniquely to an involution $\sigma \in \mathfrak{S}_j$.  Note that $k+1$ and $(k+1)'$  are not in any marked  blocks.  Now identifying $k+1$ with $(k+1)'$, merging the block $X_i$ with $X'_{\sigma(i)}$
    by adding an edge between their maximal elements, and then re-labeling the vertices
    as $-k, \dots, -1, 0, 1, \dots, k$ from left to right, 
    we get a symmetric partition of the set $[-k, k]$.  The blocks containing both positive and negative integers but not $0$ come from the marked blocks. 
\end{proof} 
\begin{example}
    This example illustrates the construction in the proof of Theorem \ref{thm: number of odd length SVT}. Let $k = 7$.
    Take $\bsy{B}^*=\{ 2^*\ | \ 1,3 \ |\ 6^*  \ | \ 5,7^*  \ | \ 4,8\ \} $ and $T= \begin{ytableau}
1 & 3 \\
2
\end{ytableau}.$
The marked blocks are $X_1=\{2\}$, $X_2=\{6\}$ and $X_3=\{5,7\}$. 
The involution determined by $T$ is $\sigma=213$. 
Hence we merge the blocks $X_1$ with $X_2^\prime$, 
$X_2$ with $X_1^\prime$, and $X_3$ with $X_3^\prime$, and identify the element $8$ of $\bsy{B}^*$ with the element $8^\prime$, as shown in the figure below. 
Finally relabeling the vertices from left to right as  $-7, \dots, -1, 0, 1, \dots, 7$,  we get the symmetric partition $\{ -7,-5\mid  -6,2\mid -4,0,4\mid -3,-1,1,3 \mid -2,6 \mid 5,7\} $ of the set $[-7,7]$. 
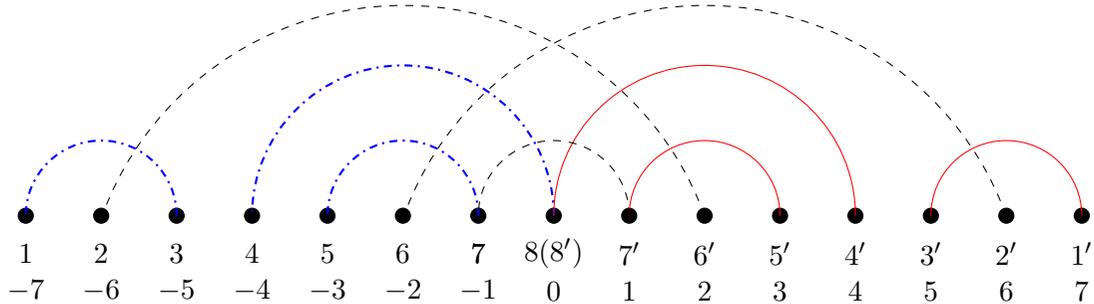
\begin{center}
    \begin{tikzpicture}
     \foreach \r in {1, ..., 15}{
      \filldraw (\r,0) circle (0.1); 
      }
      \foreach \r in {1, ..., 7}{ 
      \node at (\r, -0.5) {$\r$};
      \node at (16-\r, -0.5) {$\r^\prime$}; 
      }
      \foreach \r in {-7, ..., 7}{
       \node at (\r+8, -1) {$\r$}; 
      }
       \node at (8, -0.5) {$8(8^\prime)$}; 
       \node at (7, -0.5) {$7$}; 
      \draw[color=blue,thick,dash dot] (3,0) arc (0:180:1);  
      \draw[color=blue,thick,dash dot] (8,0) arc (0:180:2); 
      \draw[color=blue,thick,dash dot] (7,0) arc (0:180:1); 
      \draw[red]  (15,0) arc (0:180: 1); 
      \draw[red] (12,0) arc (0:180:2); 
      \draw[red] (11,0) arc (0:180:1); 
      \draw[dashed] (10,0) arc (20:160: 4.25); 
      \draw[dashed] (14,0) arc (20:160:4.25); 
      \draw[dashed] (9,0) arc (0:180:1); 
    \end{tikzpicture}
\end{center}     
\end{example}

\vanish{ 
Sequence $(g_{k+\frac{1}{2}})_{k=0}^\infty$ is  \href{https://oeis.org/A080337}{A080337} in OEIS, with the initial values $1, 3, 12, 59, 339\dots$ 
Indeed, combining with the sequence \href{https://oeis.org/A0002872}{A002872} and considering $g_0, g_{\frac{1}{2}}, g_1, g_{1\frac{1}{2}}, g_2, \dots$ , we get the sequence  \href{https://oeis.org/A080107}{A080107} in OEIS, which is described as the number of fixed points of set partitions under the involution  $i \leftrightarrow n+1-i$. It is clearly equivalent to symmetric set partitions.  
}

\medskip
\noindent 
\textsc{Remark.} 
For two sequences $(a_k)_{k=0}^\infty $ \and $(b_k)_{k=0}^\infty$, we say that 
$(b_k)_{k=0}^\infty $ is the binomial transform of  $(a_k)_{k=0}^\infty$ if $b_k=\sum_{j=0}^k \binom{k}{j}a_j$ for all $k$. In this case, we have $\sum_k b_k x^k/k! = e^x \sum_k a_k x^k/k!$. 

As an example, we note  that the sequence $(g_{k+\frac{1}{2}})_{k=0}^\infty$ is the binomial transform of $(g_k)_{k=0}^\infty$, that is, $g_{k+\frac{1}{2}}=\sum_{j=0}^k \binom{k}{j} g_j$.  Indeed, in order to form a symmetric partition of $[-k, k]$, we need exactly one self-symmetric block that contains 0, and form a symmetric set partition  with the remaining numbers.  Suppose there are $i$ positive integers (and hence $i$ negative integers by symmetry) in the block containing 0,  then there are $\binom{k}{i}$ ways to form such a block and there are $g_{k-i}$ ways to form a symmetric partition of the set of the remaining numbers. Hence  $ g_{k+\frac{1}{2}}=\sum_{i=0}^k \binom{k}{i} g_{k-i} = \sum_{j=0}^k\binom{k}{j} g_{j}$. A similar relation holds for several other pairs from the enumeration of \vtx. See Section \ref{sec:final} for a summary.

\subsection{The sum \texorpdfstring{$\sum_\mu g_{k+\frac{1}{2}}(\mu) f^\mu$}{}}  \mbox{} 

Let $u_{k+\frac{1}{2}} =\sum_\mu g_{k+\frac{1}{2}}(\mu) f^\mu$. 
Similar to the discussion for $u_k$ in Theorem \ref{thm:pk}, one can prove that  $u_{k+\frac{1}{2}}$ counts the number of partly ordered set partitions of $[k+1]$, with the additional condition that the block containing $k+1$ is not marked and hence  not in the linear order.

The initial values of $u_{k+\frac{1}{2}}$ are $1, 3, 12, 61, 381, 2854, \dots$. This sequence is not included in OEIS, but it is easy to obtain the following properties. 
\begin{enumerate}[(a)]
\item 
$\displaystyle u_{k+\frac{1}{2}} = \sum_{j \geq0}  j! \tilde B(k+1,j)$, 
where $\tilde B(k+1,j)$ is given in \eqref{B'_kj}. 

\item 
The sequence $(u_{k+\frac{1}{2}})_{k=0}^\infty$ is the binomial transform of the sequence $(u_k)_{k=0}^\infty$  (\href{https://oeis.org/A059099}{A059099}). 
 The binomial coefficient in the binomial transform corresponds to the number of ways to form the block containing the largest element $k + 1$.
 The exponential generating function of 
$(u_{k+\frac{1}{2}})_{k=0}^\infty$ is 
\[
\sum_{k\geq 0} u_{k+\frac{1}{2}} \frac{x^k}{k!} 
=e^x \cdot \frac{\exp(e^x-1) }{2-e^x}. 
\]
\end{enumerate}

\subsection{Product of \texorpdfstring{$g_{k+\frac{1}{2}}(\mu)$}{} 
with Schur functions \texorpdfstring{$\sum_\mu g_{k+\frac{1}{2}}(\mu) s_\mu(\bsy{x})$}{}} \mbox{} 

\label{ssec:schur2'} 

The following theorem is an analog of Theorem \ref{Thm:gk-schur1}.

\begin{theorem}   
Let  $k \geq 0$ be an integer.  
We have the  identity
    \begin{equation} \label{gk'-schur1}
         \sum_\mu g_{k+\frac{1}{2}}(\mu) s_\mu(\bsy{x}) = \sum_{j \geq 0} \tilde B(k+1, j) h_1(\bsy{x})^j,
    \end{equation}
    where $\tilde B(k,j)$ is  the number of set partitions of $[k+1]$ with $j$ marked blocks such that the block containing $k+1$ is not marked. 
\end{theorem}
We omit the proof because it is almost identical to that of Theorem \ref{Thm:gk-schur1}.

The principal specialization of Identity \eqref{gk'-schur1} is
\begin{eqnarray*}
    \sum_\mu g_{k+\frac{1}{2}} s_\mu(1, q, \dots, q^m) 
    &=& \sum_{j \geq 0} \tilde B(k+1, j) [m]_q ^j, 
\end{eqnarray*}

\subsection{Product of \texorpdfstring{$m_{n, k-\frac{1}{2}}^{\lambda}$}{} 
with Schur functions \texorpdfstring{$\sum_\lambda m_{k-\frac{1}{2}}(\mu) s_\lambda(\bsy{x})$}{}} \mbox{} 
\label{ssec:schur2''} 

  The following theorem is an analog of Theorem  \ref{thm:gk-schur2}.
    \begin{theorem} \label{thm:gk'-schur2}
    For all  positive integers $n$ and $k$, 
    \begin{eqnarray} \label{gk'-schur2}
        \sum_{\lambda \in \Lambda_{n-1}^k}  m_{n,k-\frac{1}{2}}^\lambda  s_\lambda(\bsy{x})
        &=& \sum_{j=1}^{\min(n,k)} S(k,j) h_{(n-j, 1^{j-1})}  (\bsy{x}).  
    \end{eqnarray}
\end{theorem}
\begin{proof}
  The proof is  similar to that of Theorem \ref{thm:gk-schur2} except that we need to explain what one gets when applying the bijection constructed in \cite{BHH17} to 
 $n$-vacillating tableaux of shape $\lambda \in \Lambda_{n-1}^k$ and length $2k-1$. 
 Analyzing the bijection of \cite{BHH17}, we see that in this case  we would obtain a
 pair $(\bsy{B}, \hat{T})$ where $\bsy{B}$ is a partition of $[k]$ into $j$ blocks   for some $j \leq n$, and $\hat{T}$ is a semistandard tableau of shape $\lambda \vdash (n-1)$ filled with $n-j$ zeros and $j-1$ distinct positive integers, which are  the maximal elements of the blocks of $\bsy{B}$ except $k$. 

 Using the same argument as for Theorem \ref{thm:gk-schur2}, we can represent an $n$-vacillating tableaux of shape $\lambda \in \Lambda_{n-1}^k$ and length $2k-1$ by a pair $(\bsy{B}, T)$ where $\bsy{B} \in \Pi(k,j)$, 
 $T \in \mathcal{SYT}(\lambda)$, and $\lambda_1 \geq n-j$.
 (Note that  $\lambda \vdash (n-1)$ instead of $\lambda \vdash n$.) 
 The remaining of the proof is the same as that of Theorem~\ref{thm:gk-schur2}. 
\end{proof}

The principal specialization of Identity \eqref{gk'-schur2} is
\begin{eqnarray*}
    \sum_{\lambda \in \Lambda_{n-1}^k}  m_{n,k-\frac{1}{2}}^\lambda  s_\lambda(1, q, \dots, q^m)
        &=& \sum_{j=1}^{\min(n,k)} S(k,j)
        \qbinom{m+n-j-1}{n-j} [m]_q^{j-1}. 
\end{eqnarray*}



\section{Limiting vacillating tableaux  of even length } \label{sec:lvt1}  

The following two sections focus on  limiting \vtx. 
Restricting the map $\psi$ of Subsection \ref{ssec:vt}  to limiting vacillating tableaux, we proved in \cite[Prop.7]{BHPYZ} that
the map $\psi$ induces a bijection 
\[
\mathcal{SVT}_k(\mu) \longleftrightarrow \Pi(k, |\mu|) \times 
\mathcal{SYT}(\mu), 
\]
 where $\Pi(k, j)$ is the set of partitions of $[k]$ with exactly $j$ blocks and 
 $\mathcal{SVT}(\mu)$ is the set of SYTs of shape $\mu$. 
  Indeed, the set partition $\bsy{B}$ and the tableaux $T_k$ constructed when applying the algorithm of $\psi$ to a limiting \vt \ satisfy the property that content$(T_k)=\max(\bsy{B})$. 
 Hence $\bsy{B}^*$ would have all its blocks marked. Therefore 
 we can represent a  limiting vacillating tableau of shape $\mu$ and length $2k$ simply by a pair 
$(\bsy{B}, T)$, where $\bsy{B}$ is a set partition of $[k]$ with $j=|\mu|$ blocks  and $T$ is a SYT of shape $\mu$.

\begin{example} Let $k = 4$ and $\mu = (2, 1)$. There are twelve limiting vacillating tableaux of shape $ (2,1)$ and length $8$. The corresponding pairs $(\bsy{B}, T)$ are given in Table~\ref{tab:1}.

\renewcommand{\arraystretch}{1.5} 
\begin{table}[htp]
\ytableausetup{smalltableaux} 
    \begin{tabular}{|c|c|} \hline 
   limiting vacillating tableau   &  $ (\bsy{B}, T) $  \\ \hline 
     $(\emptyset\,, \emptyset\,,  \ydiagram{1}\,,\emptyset\,, \ydiagram{1}\,, \ydiagram{1}\,, \ydiagram{2}\,, \ydiagram{2}\,, \ydiagram{2,1}\ )$ & $(\{1,2 \mid 3 \mid 4\}, \ytableaushort{12,3}\ )$ \\ [7pt]\hline
      $(\emptyset\,, \emptyset\,,  \ydiagram{1}\,, \ydiagram{1}\,, \ydiagram{1,1}\,, \ydiagram{1}\,, \ydiagram{2}\,, \ydiagram{2}\,, \ydiagram{2,1}\ )$& $(\{1,3 \mid 2 \mid 4\}, \ytableaushort{12,3}\ )$ \\ [7pt]\hline
     $(\emptyset\,, \emptyset\,,  \ydiagram{1}\,, \ydiagram{1}\,, \ydiagram{1,1}\,,  \ydiagram{1,1}\,, \ydiagram{2,1}\,, \ydiagram{2}\,, \ydiagram{2,1}\ )$& $(\{1,4 \mid 2 \mid 3\}, \ytableaushort{12,3}\ ) $ \\ [7pt]\hline
      $(\emptyset\,, \emptyset\,,  \ydiagram{1}\,, \ydiagram{1}\,, \ydiagram{2}\,, \ydiagram{1}\,, \ydiagram{2}\,, \ydiagram{2}\,, \ydiagram{2,1}\ )$& $(\{2,3 \mid 1 \mid 4\}, \ytableaushort{12,3}\ ) $ \\ [7pt]\hline
     $(\emptyset\,, \emptyset\,,  \ydiagram{1}\,, \ydiagram{1}\,, \ydiagram{2}\,, \ydiagram{2}\,, \ydiagram{2,1}\,, \ydiagram{2}\,, \ydiagram{2,1}\ )$& $(\{2,4 \mid 1 \mid 3\}, \ytableaushort{12,3} \ )$ \\ [7pt]\hline
     $(\emptyset\,, \emptyset\,,  \ydiagram{1}\,, \ydiagram{1}\,, \ydiagram{2}\,, \ydiagram{2}\,, \ydiagram{3}\,, \ydiagram{2}\,, \ydiagram{2,1}\ ) $& $(\{3,4 \mid 1 \mid 2\}, \ytableaushort{12,3}\ )$ \\ [7pt]\hline
    $( \emptyset\,, \emptyset\,,  \ydiagram{1}\,,\emptyset\,, \ydiagram{1}\,, \ydiagram{1}\,, \ydiagram{1,1}\,, \ydiagram{1,1}\,, \ydiagram{2,1} \ ) $      & $(\{1,2 \mid 3 \mid 4\}, \ytableaushort{13,2}\ )$ \\ [7pt]\hline
     $( \emptyset\,, \emptyset\,,  \ydiagram{1}\,, \ydiagram{1}\,, \ydiagram{1,1}\,, \ydiagram{1}\,, \ydiagram{1,1}\,, \ydiagram{1,1}\,, \ydiagram{2,1}\ ) $ & $(\{1,3 \mid 2 \mid 4\}, \ytableaushort{13,2}\ )$ \\ [7pt]\hline
     $( \emptyset\,, \emptyset\,,  \ydiagram{1}\,, \ydiagram{1}\,, \ydiagram{1,1}\,,  \ydiagram{1,1}\,, \ydiagram{1,1,1}\,, \ydiagram{1,1}\,, \ydiagram{2,1}\ ) $& $(\{1,4 \mid 2 \mid 3\}, \ytableaushort{13,2}\ )$ \\ [7pt]\hline
     $( \emptyset\,, \emptyset\,,  \ydiagram{1}\,, \ydiagram{1}\,, \ydiagram{2}\,, \ydiagram{1}\,, \ydiagram{1,1}\,, \ydiagram{1,1}\,, \ydiagram{2,1}\ )$& $(\{2,3 \mid 1 \mid 4\}, \ytableaushort{13,2}\ )$ \\ [7pt]\hline
     $( \emptyset\,, \emptyset\,,  \ydiagram{1}\,, \ydiagram{1}\,, \ydiagram{2}\,, \ydiagram{2}\,, \ydiagram{2,1}\,, \ydiagram{1,1}\,, \ydiagram{2,1}\ )$& $(\{2,4 \mid 1 \mid 3\}, \ytableaushort{13,2}\ )$ \\ [7pt]\hline
    $( \emptyset\,, \emptyset\,,  \ydiagram{1}\,, \ydiagram{1}\,, \ydiagram{1,1}\,, \ydiagram{1,1}\,, \ydiagram{2,1}\,, \ydiagram{1,1}\,, \ydiagram{2,1}\ ) $ & $(\{3,4 \mid 1 \mid 2\}, \ytableaushort{13,2}\ )$ \\ [7pt]\hline
    \end{tabular} 
        \caption{The limiting vacillating tableaux of shape $ (2,1)$ and length $8$.}
    \label{tab:1}
\end{table}

\end{example}

\subsection{Formula for \texorpdfstring{$a_k(\mu)$}{}} \mbox{} 

Let $a_k(\mu)$ be the number of limiting vacillating tableau of shape $\mu$ and length $2k$. It follows from the above bijection that
for $k \geq 0$, 
\begin{eqnarray} \label{a_k-mu}
    a_k(\mu) &=& S(k,|\mu|) f^{\mu},
\end{eqnarray}
where $S(k,j)$ is the Stirling number of the second kind and $f^\mu$ is the number of SYT of shape $\mu$. 
Note that $a_0(\emptyset)=1$ and  $a_k(\emptyset)=0$ for $k >0$.

\subsection{The sum \texorpdfstring{$\sum_\mu a_{k_1}(\mu) a_{k_2}(\mu)$ }{}}  
\mbox{} 

Let $k_1, k_2$ be two positive integers and $k=k_1+k_2$. We say that a set partition  of $[k]$ is 
$(k_1, k_2)$-connecting if for any block $X$ of the partition, 
$\min{X} \leq k_1 < \max{X}$, or equivalently,  
$X \cap [k_1] \neq \emptyset$ and 
$X \cap \{k_1+1, \dots, k_1+k_2\} \neq \emptyset$. 

\begin{example}
    Let $k_1=k_2=2$. Then the set partitions $\{1,3\ |  \ 2, 4\}$  
    and  $\{1,4 \ |  \ 2, 3\} $ are $(2,2)$-connecting 
    while  $\{1, 2\ |  \ 3,4\}$ is not.  There are altogether three $(2,2)$-connecting set partitions of 
    $[k]$ for $k=4$; the third one is $\{ 1,2,3,4\}$.  
\end{example}

\begin{theorem} \label{connect-sp}
For all integers $k_1, k_2 \geq 1$, 
    the sum $\sum_\mu a_{k_1}(\mu) a_{k_2}(\mu)$  is the number of $(k_1,k_2)$-connecting set partitions of $[k_1 +k_2]$. 
\end{theorem}
\begin{proof}
    For $i=1, 2$, let $P_{i}$ be a limiting vacillating tableau of shape $\mu$ and length $2k_i$. Then 
    $P_i$ can be represented as a pair $(\bsy{B}_i, T_i)$, where $\bsy{B}_i$ is a set partition of $[k_i]$  with $j=|\mu|$ blocks, and $T_i$ is a SYT of shape $\mu$. 
     Draw the standard diagram of $\bsy{B}_1$ with vertices $1, 2, \dots, k_1$, then followed by the reverse of the standard diagram of $\bsy{B}_2$ with vertices labeled $k_2', \dots, 2',1'$. 
    Assume the blocks of $\bsy{B}_i$ are $X^i_1, X^i_2, \dots, X^i_j$, ordered by  their maximal elements.    Using the RSK algorithm, 
    the pair $(T_1, T_2)$ determines uniquely  permutation $\sigma$.  Now merge the blocks of $\bsy{B}_1$ with 
    those of $\bsy{B}_2$ by connecting the maximal element of the blocks $X^1_t$ with that of  $X^2_{\sigma(t)}$, for all $t \leq j$.  Finally re-labeling the vertices $k_2', \dots, 2', 1'$ with $k_1+1, \dots, k_1+k_2$, we 
    get the $(k_1, k_2)$-connecting set partition. 

    Conversely, given a $(k_1, k_2)$-connecting set partition, one can easily recover $\bsy{B}_1, \bsy{B}_2$ and the permutation $\sigma$, and hence the SYT $T_1, T_2$. Therefore the above procedure is a bijection.   
\end{proof}

The sum $\sum_\mu a_{k_1}(\mu) a_{k_2}(\mu)$  can be computed as
\begin{equation}
    \sum_\mu a_{k_1}(\mu) a_{k_2}(\mu) = \sum_{j\geq 0}   j! S(k_1, j) S(k_2,j). 
\end{equation}
When $k_1=k_2$,  the sum is sequence \href{https://oeis.org/A023997}{A023997} in OEIS.

\subsection{The sum \texorpdfstring{$a_k=\sum_\mu a_k(\mu)$}{}}  \mbox{}
\label{ssec:ak} 

Let $a_k$ be the  number of limiting \vtx \ of length $2k$.   
From Formula \eqref{a_k-mu}, we have the following identity.
\begin{theorem} \cite{BHPYZ} \label{thm:lvt-ak}
 For $k \geq 0$, 
    \begin{eqnarray}
        a_k &= &\sum_{j=0}^{k}  \left( S(k,j)  \sum_{\mu \vdash j} f^{\mu} \right) = \sum_{j=0}^{k} S(k,j)  I_j,    
    \end{eqnarray}
   where $I_j$ is the number of involutions in $\mathfrak{S}_j$ and $I_0=1$. 
\end{theorem}
The sequence $(a_k)_{k=0}^\infty$ is  \href{https://oeis.org/A004211}{A004211}  in OEIS
with the initial values $1,1,3,11,49,257 \dots$. 
It has the following combinatorial interpretations relating to set partitions.  For all integers $k \geq 1$, 
 \begin{enumerate}[(a)]
     \item $a_k$ is the number of set partitions of  $[k]$ such that each element is colored either red or blue, and for each block the minimal element is colored red, \cite{GP11}. 
 \item $a_k$ is the number of pairs $(\bsy{B}, \sigma)$ where $\bsy{B}$ is a set partition of $[k]$ and
 $\sigma$ is an involution defined on the  blocks of $\bsy{B}$, \cite{BHPYZ}. 
 \end{enumerate}

The proof of Theorem \ref{connect-sp} provides a new combinatorial interpretation. In that proof  when 
$(\bsy{B}_1, T_1)= (\bsy{B}_2, T_2)$, the pair $(T_1, T_2)$ corresponds to an involution.  Hence the construction  yields a symmetric set partition of $[2k]$, (symmetric under the map $i \leftrightarrow 2k+1-i$). 

\begin{corollary}
    For  $k \geq 1$, the sum $a_k=\sum_\mu a_k(\mu)$ is the number of symmetric $(k,k)$-connecting set partitions of $[2k]$. 
\end{corollary}

\subsection{The sum \texorpdfstring{$\sum_\mu a_k(\mu) f^\mu$}{}} \mbox{} 

\begin{theorem} \label{thm: Fubini}
    Let  $v_k=\sum_\mu a_k(\mu) f^\mu$. Then $v_k$ is the $k$-th Fubini
    number  that counts the number of ordered set partitions of $[k]$. 
\end{theorem}
\begin{proof}
    This is because $v_k$ counts  the number of   triples $(\bsy{B}, T, S)$ where $\bsy{B}$ is a partition of $[k]$ with 
    $j=|\mu|$ blocks, and $T$ and $S$ are SYT of shape $\mu$. 
    Via the RSK algorithm, the pair $(T, S)$ corresponds to  a permutation in $\mathfrak{S}_j$, which gives the ordering on the blocks of $\bsy{B}$. 
\end{proof}


\subsection{Product of \texorpdfstring{$a_k(\mu)$}{} with Schur functions \texorpdfstring{$\sum_\mu a_k(\mu) s_\mu(\bsy{x})$}{}} \mbox{}

Next, we connect the limiting vacillating tableaux to symmetric functions. 

\begin{theorem} \label{lvt:schur}
For all integers $k \geq 0$,  
    \begin{eqnarray} \label{ak-schur}
        \sum_\mu a_k(\mu) s_\mu(\bsy{x})  = \sum_{j \geq 0}  S(k,j) h_1(\bsy{x})^j.
    \end{eqnarray}
\end{theorem}
\begin{proof}
    The left side of Identity \eqref{ak-schur} is a summation of $\bsy{x}^S$ over the set of triples $(\bsy{B}, T, S)$ where $\bsy{B}$ is a set partition of $[k]$ into $j$ blocks, $T$ is a SYT of some shape $\mu$ of size $j$, and $S $ is a SSYT of shape $\mu$ filled with positive integers.   The pair $(T, S)$, under the inverse of  RSK algorithm, corresponds uniquely to  sequence 
    of positive integers  $\bsy{i}=(i_1, i_2, \dots, i_j)$, where 
    \[
    \bsy{x}^S = x_{i_1} x_{i_2} \cdots x_{i_j}. 
    \]
    Hence 
    \[
    \sum_{(\bsy{B}, T, S)} \bsy{x}^S = \sum_j S(k, j) 
    \sum_{\bsy{i} \in (\mathbb{Z}^+)^j} x_{i_1} x_{i_2} \cdots x_{i_j} 
    = \sum_j S(k,j) \left(\sum_i x_i\right)^j 
    = \sum_j S(k,j) h_1(\bsy{x})^j.
    \qedhere\]    
\end{proof}

The principle specialization of Identity \eqref{ak-schur} is  
\begin{equation}  \label{ak-Schur-q}
\sum_\mu a_k(\mu) s_\mu(1, q, \dots, q^{n-1}) 
=\sum_{j\geq 0 } S(k,j) [n]_q^j. 
\end{equation}

In particular, for $q=1$  we get 
\begin{eqnarray} \label{schur-q=1}
\sum_\mu a_k(\mu) s_\mu(1^n)= \sum_{j \geq 0}  S(k,j) n^j, 
\end{eqnarray} 
where 
\[
s_\mu(1^n) = \prod_{u \in \mu} \frac{n+c(u)}{h(u)}. 
\]
Note that the right side of Identity \eqref{schur-q=1} is the evaluation of the $k$-th Bell polynomial $B_k(x):= \sum_j S(k,j) x^j$ at $x=n$. 
It is interesting to compare 
Identities \eqref{ak-Schur-q}, \eqref{schur-q=1}, and
the $q$-hook-content formula 
\eqref{q-hook-content}  with the following identity of Halverson and Thiem 
\cite[Corollary 2.4]{HT10}: for all positive integers $n, k$, 
\begin{eqnarray} \label{ht-cor2.4} 
\sum_{\lambda \in \Lambda_n^k} f^\lambda(q) m_{n,k}^\lambda = \sum_{j} S(k, j) [n]_q [n-1]_q \cdots [n-j+1]_q, 
\end{eqnarray} 
where $f^\lambda(q)$ is a $q$-analog of $f^\lambda$ and is given by the $q$-hook length formula
\[
f^\lambda(q) = q^{b(\lambda)} [n]_q! \prod_{ u \in \lambda} \frac{1}{[h(u)]_q}. 
\]
When $q=1$, \eqref{ht-cor2.4} reduces to Identity \eqref{n-vt} by the equation $x^k=\sum_{r} S(k,r) (x)_r$, where $(x)_r=x(x-1) \cdots (x-r+1).$

\vskip 1em 
\section{Limiting vacillating tableaux  of odd length}
  \label{sec:lvt2} 

When applying the map $\psi$ to a limiting vacillating tableau of length $2k$,
if we stop at the step $2k-1$ with the index $k-\frac{1}{2}$ but skip the last step, then we will have the edge set $E_{k-\frac{1}{2}}=E_k$, but $T_{k-\frac{1}{2}}$ is a SYT whose content does not contain $k$. Adjusting index we get that limiting vacillating tableaux of 
shape $\mu$ and length $2k+1$  are in one-to-one correspondence with the pairs $(\bsy{B}, T)$, where $\bsy{B}$ is a set partition of $[k+1]$ with $j=|\mu|+1$ blocks, and $T$ is a SYT of  shape $\mu$.

\subsection{Formula for \texorpdfstring{$a_{k+\frac{1}{2}}(\mu)$}{}}  \mbox{} 

Let $a_{k +\frac{1}{2}}(\mu)$ be the number of limiting vacillating tableaux of shape $\mu$ and length $2k+1$. 
From the above correspondence, we have that for $k \geq 0$, 
\begin{equation} \label{for:a-half}
    a_{k+\frac{1}{2}}(\mu)= S(k+1, |\mu|+1) f^\mu. 
\end{equation}
Note that $a_{k+\frac12}(\emptyset)=1$. 

\subsection{The sum \texorpdfstring{$\sum_\mu a_{k_1+\frac{1}{2}}(\mu) a_{k_2+\frac{1}{2}}(\mu)$ }{}}  \mbox{} 

Let  $\ell \leq k$ be positive integers. We say that a set partition  of $[k]$ is 
$\ell$-connecting if for any block $X$ of the partition, 
$\min(X) \leq \ell  \leq \max(X)$. Equivalently, a set partition of $[k]$ is $\ell$-connecting if and only if for any block $X$,
$X \cap [\ell] \neq \emptyset$ and $X \cap \{ \ell, \dots, k\} \neq \emptyset$. 

\begin{theorem} For all integers $k_1, k_2 \geq 0$, 
    the sum $\sum_\mu a_{k_1+\frac{1}{2}}(\mu) a_{k_2+\frac{1}{2}}(\mu)$ is the number of $(k_1+1)$-connecting set partitions of $[k_1+k_2+1]$.   
\end{theorem}
\begin{proof}
The proof is similar to that of Theorem \ref{connect-sp}, except 
that we start with two set partitions of $[k_1+1]$ and $[k_2+1]$ with $j+1$ blocks,  respectively, where $j=|\mu|$. 
 Note that the block containing $(k_i+1)$ is the $(j+1)$-th block of $\bsy{B}_i$ for $i = 1, 2$. 
After merging the blocks according to the permutation determined by the two SYTs and identifying the element $k_1+1$ of $\bsy{B}_1$ with $k_2'+1$ of $\bsy{B}_2$, we get a set partition of $[k_1+k_2+1]$.
\end{proof}

The sum $\sum_\mu a_{k_1+\frac{1}{2}}(\mu) a_{k_2+\frac{1}{2}}(\mu)$  can be computed as 
\begin{equation}
    \sum_\mu a_{k_1+\frac12}(\mu) a_{k_2+\frac12}(\mu) = \sum_{j \geq 0}   j! S(k_1+1, j+1) S(k_2+1,j+1). 
\end{equation}
When $k_1=k_2$,  the sum  is the sequence \href{https://oeis.org/A014235}{A014235} in OEIS. 

With a similar argument, we also have the following results. 
\begin{corollary}

\begin{enumerate}[(i)]
     \item  For all integers $k_1, k_2 \geq 0$, 
     $ \sum_\mu a_{k_1+\frac{1}{2}}(\mu) a_{k_2}(\mu)$ 
     is the number of set partitions of $[k_1+k_2+1]$ such that 
      $\min(X) \leq k_1+1 \leq \max(X)$ for all blocks $X$, and
      $k_1+1$ is the maximal element in its block. This sum can be computed as
     \[ 
     \sum_\mu a_{k_1+\frac{1}{2}}(\mu) a_{k_2}(\mu)
      = \sum_{j \geq 0}  j! S(k_1+1, j+1)S(k_2, j). 
     \]

    \item   For all integers $k_1 \geq 1$ and $k_2 \geq 0$, 
    $ \sum_\mu g_{k_1}(\mu)a_{k_2}(\mu)
    = \sum_{\mu} g_{k_1-\frac{1}{2}}(\mu) a_{k_2+\frac{1}{2}}(\mu)$  is the number of set partitions of $[k_1+k_2]$ such that  $\min(X) \leq k_1$ for all blocks $X$. This sum can be computed as 
    \begin{eqnarray*} 
     \sum_\mu g_{k_1}(\mu)a_{k_2}(\mu)
    &=& \sum_{\mu} g_{k_1-\frac{1}{2}}(\mu) a_{k_2+\frac{1}{2}}(\mu) \\ 
    &=& \sum_{j \geq 0}  j! B(k_1, j) S(k_2,j) \\
    &=&  \sum_{j \geq 0}   j! \tilde B(k_1, j) S(k_2+1, j+1). 
    \end{eqnarray*}

    \item For all integers $k_1, k_2 \geq 0$,  $ \sum_\mu g_{k_1+\frac{1}{2}}(\mu) a_{k_2}(\mu) $
     is the number of  set partitions of $[k_1+k_2+1]$ such that  $\min(X) \leq k_1+1$ for all blocks $X$, and $k_1+1$ is the maximal element in its block. This sum can be computed as 
    \[
    \sum_\mu g_{k_1+\frac{1}{2}}(\mu) a_{k_2}(\mu) 
    =\sum_{j\geq 0} j! \tilde B(k_1+1, j) S(k_2, j). 
    \]

    \item   For all integers $k_1, k_2 \geq 0$,   $\sum_{\mu} g_{k_1}(\mu) a_{k_2+\frac{1}{2}}(\mu) $ is the number of set partitions of $[k_1+k_2+1]$  such that $\min(X) \leq k_1+1$ for all blocks $X$, and $k_1+1$ is the minimal element in its block. This sum can be computed as 
    \[
    \sum_{\mu} g_{k_1}(\mu) a_{k_2+\frac{1}{2}}(\mu)
    =\sum_{j \geq 0}  j! B(k, j) S(k+1, j+1). 
    \]
\end{enumerate}
\end{corollary}

\subsection{The sum \texorpdfstring{$\sum_\mu a_{k+\frac{1}{2}}(\mu)$}{}}  \mbox{}

Let $a_{k+\frac{1}{2}}= \sum_\mu a_{k+\frac{1}{2}}(\mu)$ be the number of limiting \vtx\ of shape $\mu$ and length $2k+1$. 
The next theorem  follows immediately from Formula \eqref{for:a-half}.  
\begin{theorem} \label{sum:a_half}
For $k \geq 0$,
    \begin{equation} \label{eq:a_half}
        a_{k+\frac{1}{2}} = \sum_{j=0}^k S(k+1, j+1) \sum_{\mu \vdash j} f^\mu = \sum_{j=0}^{k} S(k+1,j+1) I_j
    \end{equation}
    where $I_j$ is the number of involutions in $\mathfrak{S}_j$ and $I_0=1$. 
\end{theorem}

The numbers of limiting vacillating tableaux of  odd and even 
length, respectively, are closely related. 
\begin{theorem} \label{thm:a_half}
   The sequence $(a_{k+\frac{1}{2}})_{k=0}^\infty$  is the binomial transform of $(a_k)_{k=0}^\infty$. That is, 
   \begin{equation}
       a_{k+\frac{1}{2}} = \sum_{r=0}^k \binom{k}{r} a_r. 
   \end{equation}
\end{theorem}
\begin{proof}
  $a_{k+\frac12}$ counts the number of pairs of  $(\bsy{B}, T)$, where $\bsy{B}$ is a set partition of $[k+1]$ with $j+1$ blocks, and $T$ is a SYT of  shape $\mu$ with $\mu \vdash j$. 
    For a partition of $[k+1]$ with $j+1$ blocks, assume the block containing $k+1$ has $k+1-r$ elements. There are $\binom{k}{r}$ ways to pick those elements. The remaining $r$ elements then form a partition with  $j$ blocks, which, together with the SYT $T$ of shape  $\mu$, correspond to a limiting vacillating tableau of  length $2r$. 
\end{proof}

The sequence $(a_{k+\frac{1}{2}})_{k=0}^\infty$  is  \href{https://oeis.org/A007405}{A007405} in OEIS,
whose initial terms are $1, 2, 6, 24, 116, 648, \dots$. 
Interestingly, $(a_k)_{k=0}^\infty$ can  also be obtained from  the binomial transform of $(a_{k+\frac{1}{2}})_{k=0}^\infty$.
\begin{theorem} \label{thm:a_half2}
For $k \geq 0$,
    \begin{equation}
        a_{k+1} = \sum_{j=0}^k \binom{k}{j} a_{j+\frac{1}{2}}. 
    \end{equation}
\end{theorem}
\begin{proof}
    By Equation \eqref{eq:a_half} a combinatorial interpretation of  $a_{j+\frac{1}{2}}$ is the number of pairs $(\bsy{B}', \sigma')$, where $\bsy{B}'$ is a set partition of $[j+1]$ with $r$ blocks and $\sigma' \in \mathfrak{S}_r$ is 
    an involution on the blocks of $\bsy{B}'$ such that the block containing $j+1$ is a 1-cycle.
    
    Let $(\bsy{B}, T)$ be a pair counted by $a_{k+1}$, that is, $\bsy{B}$ is a partition of $[k+1]$ whose blocks 
    are $X_1, X_2, \dots, X_t$ (ordered by their maximal elements), and 
    $T$ is a SYT of some shape $\lambda$  with $|\lambda|=t$.  
    Let $\sigma$ be the involution corresponding to $(T, T)$ under the RSK algorithm. 
    Note that $X_t$ is the block of $\bsy{B}$ containing the element $k+1$. Assume $|X_t|=1+k-j$.  Remove the block $X_t$ from $\bsy{B}$ and introduce a new element $\star$. 
    Using the following steps, we construct a set partition $\bsy{B}'$ on $([k+1]-X_t) \cup \{\star\}$ together with an involution     $\sigma'$ on the blocks of $\bsy{B}'$, such that the block containing $\star$ is a 1-cycle of $\sigma'$. 

    \begin{enumerate}[(i)]
        \item  If $\sigma(t)=t$, that is,  $X_t$ is not paired with another block of $\bsy{B}$ under $\sigma$, then   $\{ \star\} $ is the $t$-th block of $\bsy{B}'$, and this block is a 1-cycle of $\sigma'$. 
        \item If $\sigma(t)=r$ for some $r < t$, that is, $(X_t, X_r)$ is a 2-cycle  of $\sigma$, then add $\star$ to the block $X_r$ to form the $r$-th block of $\bsy{B}'$, and let $\sigma'(r)=r$.  
        \item  For any block $X_a$ of $\bsy{B}$ with $\sigma(a)=b$ where  $a, b \neq t$, 
        $X_a$ and $X_b$ are both blocks of $\bsy{B}'$ with $\sigma'(a)=b$.  Note that it is possible that $a=b$ in this step. 
    \end{enumerate}
    The resulting pairs $(\bsy{B}', \sigma')$  are counted by $a_{j+\frac{1}{2}}$.    There are $\binom{k}{j}$ ways to choose the elements in $X_t$. Summing over $j$ completes the proof. 
\end{proof}

\begin{example} Let $k = 8$.
As an example, assume that $\bsy{B} = \{2,5 \mid 1,3,6 \mid 4,8\mid 7,9\}$,  where the blocks are ordered by their maximal elements. If $\sigma = (1)(23)(4)$ in cycle notation, then the operation yields the set partition $\bsy{B}' =\{2,5 \mid 1,3,6 \mid 4,8\mid \star\}$ with involution $\sigma' = (1)(23)(4)$.
If $\sigma = (12)(34)$ in cycle notation, then the operation yields the set partition $\bsy{B}' =\{2,5 \mid 1,3,6 \mid 4,8, \star\}$ with involution $\sigma' = (12)(3)$.
\end{example}

Combining Theorems \ref{thm:a_half} and \ref{thm:a_half2}, one can derive that $(a_k)_{k=0}^\infty$ satisfies the recurrence relation 
\[
a_{k+1}= \sum_{j=0}^k \binom{k}{j} 2^{k-j}a_j. 
\]
This can be interpreted combinatorially using the model of bi-colored set partition  defined in Subsection \ref{ssec:ak}, item (a) after Theorem \ref{thm:lvt-ak}, where  $j$ is the number of elements not in the same block as $k+1$.

We remark that the sequences $(a_k)_{k=0}^\infty$ and $ (a_{k+\frac{1}{2}})_{k=0}^\infty$ are the unique pair of sequences that start at $s_0=1$ and recursively defined by 
the binomial transforms 
\[
t_k= \sum_{j=0}^k \binom{k}{j} s_j \qquad \text{ and } \qquad 
s_{k+1}= \sum_{j=0}^k \binom{k}{j} t_j. 
\]

Another combinatorial interpretation of 
$a_{k+\frac{1}{2}}$  is the number of type $B$ set partitions, which was first introduced by Reiner \cite{Reiner97} in the study of intersection lattice for the classical  reflection groups of type $B$. 

\begin{definition} \cite{Reiner97} 
A set partition of type $B$ is a partition $\pi$ of the set 
$[-k] \cup [k]$ into blocks satisfying the following conditions:
\begin{enumerate}[(a)]
\item For any block $X$ of $\pi$, its opposite $-X$ is also a block of $\pi$;
\item  There is at most one zero-block, which is defined to be a block $X$ such that $X=-X$. 
\end{enumerate} 
\end{definition}

\begin{theorem}
 For $k \geq 1$,    the integer $a_{k+\frac{1}{2}}$  is the number of type $B$ partitions
    of the set $[-k] \cup [k]$.
\end{theorem}
\begin{proof} 
    Using Equation \eqref{eq:a_half}, 
    we interpret  $a_{k+\frac{1}{2}}$ as  the number of $(\bsy{B}', \sigma')$ where $\bsy{B}'$ is a set partition of $[k+1]$ with $j+1$ blocks and $\sigma' \in \mathfrak{S}_{j}$  is  an involution of size $j$. 
    
    Take two identical copies of the set partition $\bsy{B}'$, the first on the elements $1, 2, \dots, k+1$, and the second  on the elements $-1, -2, \dots, -(k+1)$. 
    Let $X_1, \dots, X_{j+1}$ be the blocks in the first copy of $\bsy{B}'$  with $k+1 \in X_{j+1}$,  
    and $-X_1, \dots, -X_{j+1}$ the corresponding blocks in the second copy. 
    We form a type $B$ partition of $\{\pm 1, \dots, \pm k\}$ using the following steps. 
    \begin{enumerate}[(i)]
        \item If $\sigma'(a)=b$ while $a \neq b$, then merge block $X_a$ with $-X_b$, and block $X_b$ with $-X_a$. 
        \item If $\sigma'(a)=a$, then leave both $X_a$ and $-X_a$ unchanged. 
        \item  Merge $X_{j+1}$ with $-X_{j+1}$ and then remove the elements $\pm (k+1)$. 
        If the resulting set is not empty, then it is the zero block.  Otherwise, discard the empty set.         
    \end{enumerate}
    The collection of all the blocks obtained by these steps is the desired type $B$ partition. 
    We leave it to the reader to check that the above construction indeed gives a bijection.  
\end{proof}

\begin{example}
    As an example, assume  $\bsy{B}'=\{ 2 \ | \ 1,3 \ | \ 6 \ |\ 5,7\ | \ 4,8\}$, where the blocks are ordered by their maximal elements, 
    and $\sigma'=(1)(23)(4)$ in cycle notation. Then the operations yield the type $B$ partition 
    $ \{ 2\ | \ -2 \ | \ 1,3, -6\  | \ -1,-3, 6 \ | \ 5,7 \  | \ -5, -7, \ | \ 4, -4 \} $. 
\end{example}

\subsection{The sum \texorpdfstring{$\sum_\mu a_{k+\frac{1}{2}}(\mu) f^\mu$}{}}  \mbox{} 

  For a set $S$, 
a \emph{cyclically ordered set partition} of $S$ 
is a set partition $\bsy{B}$ of $S$ together with a cyclic ordering on the  blocks of $\bsy{B}$. 
For $S=[k+1]$, the number of cyclically ordered set partition can be computed by the summation 
\[
\sum_{j\geq 0}  j! S(k+1, j+1).  
\]

\begin{theorem}
    Let  $v_{k+\frac{1}{2}}= \sum_\mu a_{k+\frac{1}{2}}(\mu) f^\mu$ for $k \geq 0$. The $v_{k+\frac{1}{2}}$ is the number  of the following structures:
    \begin{enumerate}[(a)]
        \item  cyclically ordered set partitions of $[k+1]$, or 
    \item set partitions of $[k+1]$ with a linear order on all the blocks except the block containing $k+1$. 
    \end{enumerate}
     
\end{theorem}
The proof is similar to that of Theorem \ref{thm: Fubini} and hence is omitted. The sequence  
$(v_{k+\frac{1}{2}})_{k=0}^\infty$  is \href{https://oeis.org/A000629}{A000629} in OEIS, 
which is the binomial transform of  Fubini numbers.

\subsection{Product of \texorpdfstring{$a_{k+\frac12}(\mu)$}{} with Schur functions \texorpdfstring{$\sum_\mu a_{k+\frac{1}{2}}(\mu) s_\mu(\bsy{x})$}{}} \mbox{} 

The analogous result to Theorem \ref{lvt:schur} is the following identity.  Again,  we skip the proof due to its similarity to the proof of Theorem \ref{lvt:schur}.
\begin{theorem} \label{lvt:schur2}
For $k \geq 0$, 
    \begin{eqnarray} \label{ak-schur2}
        \sum_\mu a_{k+\frac{1}{2}}(\mu) s_\mu(\bsy{x})  = \sum_{j \geq 0}  S(k+1,j+1) h_1(\bsy{x})^j.
    \end{eqnarray}
\end{theorem}
Consequently, 
\[
\sum_\mu a_{k+\frac{1}{2}}(\mu) s_\mu(1, q, \dots, q^{n-1})  = \sum_{j \geq 0}  S(k+1,j+1) [n]_q^j.
\]


\section{Final remarks and Future Projects}  \label{sec:final}

Comparing  the results of Sections \ref{sec:svt1} \& \ref{sec:svt2} and those of Sections \ref{sec:lvt1} \& 
\ref{sec:lvt2}, 
we notice that  an integer sequence relating to \vtx\ 
of odd length   
is often the binomial transform of its analog 
of the \vtx \ of even length. This includes the following pairs: 
\begin{itemize}
    \item $g_k=\sum_\mu g_k(\mu)$ and $g_{k+\frac{1}{2}}=\sum_\mu g_{k+\frac{1}{2}}(\mu)$
    \item $u_k=\sum_\mu g_k(\mu) f^\mu$ and $u_{k+\frac{1}{2}}=\sum_\mu g_{k+\frac{1}{2}}(\mu) f^\mu$ 
    \item $a_k=\sum_\mu a_k(\mu)$ and $a_{k+\frac{1}{2}}=\sum_\mu a_{k+\frac{1}{2}}(\mu)$
    \item $v_k=\sum_\mu a_k(\mu) f^\mu$ and $v_{k+\frac{1}{2}}=\sum_\mu a_{k+\frac{1}{2}}(\mu) f^\mu$ 
\end{itemize}

In addition, we observe that there are several instances that a sequence relating to simplified vacillating tableaux is the binomial convolution of  the Bell numbers
and the analogous sequence relating to limiting vacillating tableaux.  The following shows such an instance. 
\begin{theorem}
We have  $  g_k=\sum_j \binom{k}{j}  B(j)a_{k-j}$,   where $B(j)$ is the $j$-th Bell number.
\end{theorem}
\begin{proof}
    A vacillating tableau of length $2k$ is represented by a pair $(\bsy{B}^*, T)$, where $\bsy{B}^*$ is a partition of $[k]$
    with $j$ marked blocks, and $T$ is a SYT of some shape $\mu$ with $|\mu|=j$. Let $S \subseteq [k]$ be the union of those un-marked blocks. Then  $\bsy{B}^*$ can be viewed as a disjoint union of two structures: 
     a set partition $\bsy{B}_1$ of $S$, and a set partition  $\bsy{B}_2$ of $[k]-S$ with exactly $j$ blocks. The number of choices for $\bsy{B}_1$ is $B(|S|)$, while the pairs $(\bsy{B}_2, T)$ correspond to limiting vacillating tableaux of length $2(k-|S|)$ that are counted by $a_{k-|S|}$. 
\end{proof}
It follows that the exponential generating function for $(g_k)_{k=0}^\infty  $ is the product of those for the Bell numbers and for $(a_k)_{k=0}^\infty$, i.e., 
\[
\sum_{k \geq 0} g_k \frac{x^k}{k!} = \exp( e^x-1) \cdot \sum_{k \geq 0} a_k \frac{x^k}{k!} 
\]
A similar relation  holds between the following pairs of sequences: 
\begin{itemize}
    \item $g_k$ and $a_k$,
    \item $g_{k+\frac{1}{2}}$ and $a_{k+\frac{1}{2}}$,
    \item $v_k$ and $u_k$, and 
    \item $v_{k+\frac{1}{2}}$ and $u_{k+\frac{1}{2}}$. 
\end{itemize}

The exponential generating functions of the above sequences are 
summarized in Table \ref{tab:gf}.  
\medskip 

\setlength{\extrarowheight}{10pt}
\begin{table}[htb]
\centering
\begin{tabular}{|p{1.8cm}|c|c|c|} 
   \hline 
   & Sequence  &  OEIS ID    &  EGF     \\ [5pt] \hline 
\multirow{4}{*}{\rotatebox[origin=c]{0}
{\shortstack{\phantom{X}\\\phantom{X}\\\phantom{X}\\\phantom{X}\\Limiting\\vacillating \\ tableaux}}}
& $a_k=\sum_\mu a_k(\mu)$  
     &  \href{https://oeis.org/A004211}{A004211}
     &  $ \displaystyle  \exp\left(\frac{e^{2x}-1}{2}\right) $ 
     \\  [10pt]\cline{2-4}
& $a_{k+\frac{1}{2}}=\sum_\mu a_{k+\frac{1}{2}}(\mu)$
     &  \href{https://oeis.org/A007405}{A007405} 
     &    $ \displaystyle 
        \exp\left( x+  \frac{e^{2x}-1}{2}\right)$ 
        \\ [10pt]\cline{2-4}
& $v_k=\sum_\mu a_k(\mu) f^\mu$ 
      & \href{https://oeis.org/A000670}{A000670} 
      &  $ \displaystyle 
  \frac{1}{2-e^x} $ 
      \\ [10pt]\cline{2-4}
& $v_{k+\frac{1}{2}}=\sum_\mu a_{k+\frac{1}{2}}(\mu) f^\mu$ 
      &  \href{https://oeis.org/A000629}{A000629} 
      &  $ \displaystyle \frac{e^x}{2-e^x}$ 
      \\ [10pt]\hline 

\multirow{4}{*}{\rotatebox[origin=c]{0}{\shortstack{\phantom{X}\\\phantom{X}\\\phantom{X}\\\phantom{X}\\Simplified\\vacillating \\ tableaux}}}
& $g_k=\sum_\mu g_k(\mu)$ 
      &  \href{https://oeis.org/A002872}{A002872} 
      & $ \displaystyle 
 \exp\left( e^x+ \frac{e^{2x}-3}{2}\right) $  
      \\ [10pt]\cline{2-4}
& $g_{k+\frac{1}{2}}=\sum_\mu g_{k+\frac{1}{2}}(\mu)$
      &  \href{https://oeis.org/A080337}{A080337} 
      & $ \displaystyle 
 \exp\left( x+ e^x+ \frac{e^{2x}-3}{2}\right) $  
      \\ [10pt]\cline{2-4}
& $u_k=\sum_\mu g_k(\mu) f^\mu$ 
     &   \href{https://oeis.org/A059099}{A059099} 
     & $ \displaystyle   \frac{\exp(e^x-1) }{2-e^x} $
     \\ [10pt]\cline{2-4}
& $u_{k+\frac{1}{2}}=\sum_\mu g_{k+\frac{1}{2}}(\mu) f^\mu$    
    &  not in OEIS
    & $ \displaystyle   \frac{\exp(x+e^x-1) }{2-e^x}$  
    \\ [10pt]\hline         
\end{tabular}
\caption{Summary of exponential generating functions.}\label{tab:gf}
\end{table}

\medskip


We conclude this paper with some final remarks and future research projects.

Using the growth diagrams, Krattenthaler \cite{Kratten06} explored the results on simplified \vtx \ and their connections to 
 crossings and nestings of set partitions  with a broader context. This context involves the enumeration of fillings of Young diagrams, where certain restrictions are imposed on the increasing and decreasing chains of the fillings. Recently, in \cite{Kratten23} Krattenthaler extended Identity \eqref{n-vt} on $n$-\vtx \ by studying the growth diagrams associated with the shape $((n+k)^n, n+k-1, \dots,n+1,  n)$.

The growth diagrams prove to be a valuable tool as they provide in-depth insights into tableau operations and facilitates visualizations of \vtx \ between arbitrary shapes $\mu$ and $\lambda$. It is intriguing to explore how the combinatorial identities considered in this paper would transform for such generalized \vtx.

 In the context of fillings of Young diagrams and other general polyominoes, the growth diagrams frequently establish sequences of integer partitions under specific restrictions, bearing resemblance to \vtx.  Notably, the shapes of these integer partitions play a significant role in characterizing the NE- and SE-chains within the fillings. This observation has been explored in various papers, including  \cite{Kratten06, Rubey11,GuoP20}. 
The paper \cite{Kratten23} leveraged this property to define a sub-family of fillings that exhibit special combinatorial structures. 
One interesting direction is to develop $q$-analogs of those combinatorial identities to 
include the statistics of NE- and SE- chains.

  Identity \eqref{n-vt} arises from the representation theory of partition algebra. 
Several other identities presented in this paper also exhibit strong indications of a connection to representation theory. Unveiling and understanding such connections would be highly valuable. It is worth noting that, although \vtx  are closely related to the partition algebra, the role of limiting \vtx \ in representation theory remains unclear and requires further exploration.

\section*{Acknowledgement}

This material is based upon work supported by the National Science Foundation under Grant No. DMS-1929284 while the authors was in residence at the Institute for Computational and Experimental Research in Mathematics (ICERM) in Providence, RI, during the Discrete Optimization: Mathematics,  Algorithms, and Computation semester program.
We thank ICERM for facilitating and supporting this research project. 

We want to thank Richard Stanley for suggesting some of the summations, which inspired this paper. We are also grateful to James Sundstrom for providing help on computer programming. 
The second author is supported by an EDGE Karen Uhlenbeck fellowship, the third author is supported by a PSC-CUNY award, and the fourth author is supported in part by the Simons Collaboration Grant for Mathematics 704276.

\end{document}